\documentclass[dvips, ,preprint, leqno]{imsart}
\RequirePackage[OT1]{fontenc} \RequirePackage{amsthm,amsmath}
\RequirePackage[numbers]{natbib}
\RequirePackage[colorlinks,citecolor=blue,urlcolor=blue]{hyperref}
\theoremstyle{plain}
\newtheorem{dfn}{Definition}[section]
\newtheorem{thm}{Theorem}[section]
\newtheorem{lemma}{Lemma}[section]
\newtheorem{cor}{Corollary}[section]

\newtheorem{rmq}{Remark}[section]
\newtheorem{epl}{Example}[section]
\newcommand{\bm}{\boldsymbol{\mu}}
\endlocaldefs
\begin{document}
\begin{frontmatter}
\title{ON A KOTZ-WISHART DISTRIBUTION:~MULTIVARIATE VARMA TRANSFORM}
\runtitle{ON A KOTZ-WISHART DISTRIBUTION }
\begin{aug}
\author{\fnms{AMADOU SARR\\Sultan Qaboos University} \ead[label=e1]{asarr7@gmail.com , asarr@squ.edu.om}}

\runauthor{AMADOU SARR}


\end{aug}

\begin{abstract}
The Wishart distribution, named after its discoverer, play a
fundamental role in classical multivariate statistical analysis. In
this paper, we introduced a generalized Wishart distribution for
which the underlying observed vectors follow a Kotz-type
distribution. Several properties of the Kotz-Wishart (KW) random
matrix and its inverted version are investigated. Explicit forms for
the probability density functions (pdf), the cumulative distribution
function (cdf), and the moment generating function (mgf) are
derived. It has also been proved that some estimating results
obtained under Efron-Morris loss function \cite{efron} are robust
under the Kotz model. Further, as a generalization of Khatri's result \cite{khatri}, the distribution function of the
smallest eigenvalue of the KW random matrix is obtained. \\On the other hand, inspired by the particular form of the
pdf of KW random matrix, we introduced a multivariate version of a
generalized Laplace transform, which is known in the literature as
Varma transform \cite{varma} (1951). Examples of M-Varma transform
are evaluated for some functions of matrix argument.
\end{abstract}

\end{frontmatter}

\section{Introduction}
Let $X_{1},\ldots,X_{n}$ be independent and identically distributed
(iid) as multivariate normal with mean $\mu$ and $p\times p$
positive definite matrix $\Sigma $, denoted as ${\cal{N}}_{p}(\mu
,\Sigma)$, $\Sigma > 0$. Let
 $\overline{X}=\frac{1}{n}\sum_{i=1}^{n}X_{i}$, with $n>p$, and the sum
of squares and products (SSP) matrix ${\bf
A}~=~\sum_{i=1}^{n}(X_{i}-\overline{X})(X_{i}-\overline{X})^{\prime}$.
Then, the random matrix ${\bf A}$ has a Wishart distribution with
$n-1$ degrees of freedom and scale parameter matrix $\Sigma$. This
is denoted by $W_{^p}(n-1, \Sigma)$.\\This matrix variate
distribution was introduced in 1928 \cite{wish} by J. Wishart, in
the context of statistical data analysis. Nowadays, the usefulness
of Wishart's finding goes beyond the multivariate statistical area.
Indeed, Wishart matrices found many applications in diverse fields
including biology, finance, physics and more recently mechanical and
electrical Engineering \cite{adhi}, and graphical models
\cite{khare} as well. Further, Wishart distributions may be viewed
as the raw material of the random matrix theory (RMT). For more
details on the matrix variate distributions, we refer to the books
by Gupta and Nagar \cite{gnagar}, Muirhead \cite{muir}, to mention a
few.\\Generally speaking, the multivariate statistical analysis,
based on the normality assumption, has a long history of
development. Several books treating this subject are available in
the statistical literature for a while \cite{ander}, \cite{muir}.
However, the normality assumption is continuously questioned for
lack of robustness. As a result, alternative models to the
multivariate normal model have received considerable attention over
the last three decades. Namely, the multivariate elliptical models
for which the observed vectors $X_{1},\ldots,X_{n}$ are identically
distributed, with a common distribution belonging to the broader
class of elliptically contoured distributions (ECD). Indeed, such a
class of distributions includes the normal distributions and the
t-distributions as special cases. Detailed studies of ECD, for both
vectors and matrices cases, may be found in \cite{fkotz},
\cite{gth}, \cite{zhang}, \cite{varga}. In our present work, we are
mainly interested in the properties of the SSP matrix, which is
generated from a multivariate Kotz-type model.
\begin{dfn}
The $p$-dimensional random vector ${\bf x}$ is said to have a
symmetric Kotz type distribution with parameters $q,\theta ,s \in
\bf{R}$, ${\mu} \in {\bf{R}}^{p}$, ${\bf \Sigma} \in {\cal
R}^{p\times p}$,~$\theta >0$, $s>0, 2q+p>2,~and~{\bf \Sigma}> 0$ if
its pdf is
\begin{equation}\label{Ko}
f({\bf x})~=~C{\mid{\bf \Sigma}\mid}^ {-\frac{1}{2}}~[({\bf
x}-{\mu})^{\prime}{\bf \Sigma}^{-1}({\bf
x}-{\mu})]^{q-1}\exp\{-{\theta}[({\bf x}-{\mu})^{\prime}{\bf
\Sigma}^{-1}({\bf x}-{\mu})]^{s}\},
\end{equation}
where \begin{equation*} C =
\frac{s{\theta}^{\frac{2q+p-2}{2s}}\Gamma(\frac{p}{2})}{{\pi}^{\frac{p}{2}}\Gamma\left(\frac{2q+p-2}{2s}\right)}
\end{equation*}
\end{dfn}
This distribution was introduced by Kotz (1975)\cite{kotz}, as a
generalization of the multivariate normal distributions. Indeed, the
particular case $q=2\theta=s=1$ in (\ref{Ko}) coincides with the
family of normal distributions. Note that this subclass of
elliptical distributions is often used when the normality assumption
is not acceptable (see Lindsey (1999)\cite{lind}). Most of the
results related to this distribution can be found in Nadarajah
(2003)\cite{nadara}. As a member of the elliptical family, this
distribution admits the stochastic representation which is
illustrated by the relation ${\bf x}~\stackrel{d}{=}~{\mu}+r{\bf
\Sigma}^{\frac{1}{2}}{\bf u}^{(p)}$, where ${\bf u}^{(p)}$ is
uniformly distributed on the surface of the unit sphere in ${\cal
R}^{p}$, and $r$ is independent of ${\bf u}^{(p)}$. Further, the
moments of $r^2$ are given by (\ref{rkot}) (see Fang et al. \cite[p.77]{fkotz}).
\begin{equation}\label{rkot}
E\left(r^{2t}\right)~=~\frac{{\theta}^{-\frac{t}{s}}\Gamma\left(\frac{2q+p+2t-2}{2s}\right)}{\Gamma\left(\frac{2q+p-2}{2s}\right)},~~~t>0.
\end{equation}
Throughout this paper, the multivariate Kotz distribution will be
denoted by $MK_{p}({\mu},{\bf \Sigma})$.
Now, we define the multivariate Kotz type model and subsequently, the Kotz-Wishart matrix.\\
Let $X_{1},\ldots,X_{n}$ be $p$-dimensional random vectors, such
that $n>p$ and ${\bf x}_{i}\sim MK_{p}({\mu},{\bf \Sigma}),$
$i=1,\ldots,n.$ Moreover, assume that ${\bf x}_{i},$ $i=1,\ldots,n$
are uncorrelated,(but not necessarily independent), and their joint
pdf is given by
\begin{equation}\label{model1}
f({\bf x}_{1},\ldots,{\bf x}_{n})~=~\frac{1}{\mid {\bf \Sigma}\mid
^{\frac{n}{2}}}h\left(\sum_{i=1}^{n}({\bf x}_{i}-{\mu})^{\prime}{\bf
\Sigma}^{-1}({\bf x}_{i}-{\mu})\right),
\end{equation}
where $h(\cdot)$ is specified by
\begin{equation}\label{modkot}
h(z)~=~\frac{s{\theta}^{\frac{2q+pn-2}{2s}}\Gamma\left(\frac{pn}{2}\right)}{\pi^{\frac{pn}{2}}\Gamma\left(\frac{2q+pn-2}{2s}\right)}z^{q-1}\exp(-{\theta}z^{s})~,~~z\geq 0.
\end{equation}
Here the parameters $q,\theta$ and $s$ are supposed to be known, and
each $p$-dimensional random vector ${\bf x}_{i}$ $(i=1,2,\cdots,n)$
has a multivariate Kotz distribution with mean vector ${\mu}$ and
covariance matrix
$\frac{{\theta}^{-\frac{1}{s}}\Gamma\left(\frac{2q+p}{2s}\right)}{p\Gamma\left(\frac{2q+p-2}{2s}\right)}{\bf
\Sigma}$. So, the relation (\ref{model1}) represents the
multivariate Kotz type model. This model has been considered by Sarr
and Gupta \cite{sarr}, in the context of estimating the precision
matrix ${\bf \Sigma}^{-1}$, under a decision-theoretic viewpoint.
More details regarding the description of the multivariate Kotz
model are available in \cite{sarr}.
\begin{dfn}{\bf Kotz-Wishart matrix}\\
The SSP matrix ${\bf
A}~=~\sum_{i=1}^{n}(X_{i}-\overline{X})(X_{i}-\overline{X})^{\prime}$,
formed from the multivariate Kotz type model, is called the
Kotz-Wishart matrix, with $n-1$ degrees of freedom and parameter
matrix ${\bf \Sigma}> 0$. This will be denoted as
\begin{equation}\label{KW}
{\bf A} \stackrel{d}{=}KW_{p}(n-1, \Sigma),
\end{equation}
with $n>p$. Its pdf will be explicitly derived later.
\end{dfn}
\begin{rmq}It is important to observe that the uncorrelated
elliptical models coincide with the usual independent model if the
underlying vectors ${\bf x}_{1},{\bf x}_{2},\ldots,{\bf x}_{n}$ are
normally distributed. Moreover, the excellent technical report
\cite{TR54} by Anderson and Fang (1982), dealing with the
uncorrelated (or dependent) elliptical models, has been a key source
of motivation for many researchers. Since then, the relevance of the
uncorrelated models has been illustrated in a paper by Kelejian and
Prucha (1985) \cite{kelej}, where the authors proved that the
uncorrelated t-model is better able to capture the heavy-tailed
behavior than an independent t-model.
\end{rmq}
The present paper is organized as follows: Section2 contains some preliminary results and notations used throughout the paper. In section3 , we evaluated some specific expected values, needed to estimate the precision matrix ${\bf \Sigma}^{-1}$ in section4. While section5 is dedicated to the derivation of the pdfs. All the results involving zonal polynomials are discussed in section6. And section7 is devoted to the definition of a generalization of the Laplace transform, namely, the M-Varma transform. Finally, some properties concerning Whittaker's functions, which are repeatedly used in our paper, are summarized in the Appendix.
\section{Preliminaries and Notations}
The preliminary results and definitions needed for the sequel are
summarized here.
\subsection{\bf Notations and Definitions}
Let ${\bf A}$ be a $p\times p$ matrix. Then its determinant is
denoted by $|\bf A|$, its transpose matrix is denoted by ${\bf
A}^{\prime}$. The trace of ${\bf A}$ is denoted by $tr {\bf
A}=a_{11}+\cdots +a_{pp}$. The norm of ${\bf A}$ is denoted by
$\|{\bf A}\|$. The exponential of the trace of ${\bf A}$ is denoted
by $etr({\bf A})=\exp(tr {\bf A})$. The notation ${\bf A}>0$ means
that ${\bf A}$ is symmetric positive definite, and ${\bf A}^{1/2}$
denotes the unique symmetric positive definite square root of ${\bf
A}$; while the notation ${\bf A}\geq 0$ means that ${\bf A}$ is
symmetric positive semi-definite. Let ${\bf M}$ and ${\bf B}$ be a
$p\times q$ and $r\times s$ matrices respectively. Then, the {\bf
Kronecker} product of ${\bf M}$ and ${\bf B}$ is denoted by ${\bf
M}\otimes{\bf B}$. We also denote by $I_{n}$ the identity matrix of
order $n$.
\begin{dfn}{\bf The multivariate gamma function}\\
The Multivariate gamma function, denoted by ${\Gamma}_{p}(\cdot)$,
is defined as:
\begin{equation}\label{mgam}
{\Gamma}_{p}(a) = \int_{A>0}etr(-A)|A|^{a-\frac{p+1}{2}}(dA),
\end{equation}
where $Re(a)>\frac{p-1}{2}$, and the integral is over the space of
positive definite $p\times p$ matrices, with respect to the Lebesgue
measure $(dA)= da_{11}da_{12}\ldots da_{pp}$. Further, the following
formula can be proved (see \cite[p.62-63]{muir}).
\begin{equation*}
\Gamma_{p}(a)~=~{\pi}^{\frac{p(p-1)}{4}}~\prod_{i=1}^{p}\Gamma(a-\frac{i-1}{2}).
\end{equation*}
\end{dfn}
When $p=1$, then (\ref{mgam}) reduces to the classical gamma
function.
\begin{dfn}{\bf Wishart and Inverted Wishart distributions}\\
A $p\times p$ symmetric random matrix ${\bf A}>0$ has a Wishart
distribution $W_{^p}(m,\Sigma)$, with m degrees of freedom, and
parameter matrix $\Sigma >0$ and $m\geq p$ if its pdf is given by:
\begin{equation}\label{Wishart}
f({\bf A})~=~\frac{{\mid{\bf
A}\mid}^{\frac{m-p-1}{2}}~etr\left(-\frac{1}{2}{\bf \Sigma}^{-1}{\bf
A}\right)}{2^{\frac{mp}{2}}{\mid {\bf
\Sigma}\mid}^{\frac{m}{2}}\Gamma_{p}(\frac{m}{2})}~,
\end{equation}
A $p\times p$ random matrix ${\bf B}>0$ is said to have the {\bf
inverted} Wishart distribution with d degrees of freedom and
positive definite $p\times p$ parameter matrix ${\bf V}$ if its pdf
is
\begin{equation}\label{invWishart}
f({\bf B})~=~\frac{{\mid{\bf
B}\mid}^{-\frac{d}{2}}~etr\left(-\frac{1}{2}{\bf V}{\bf
B}^{-1}\right)}{2^{\frac{p(d-p-1)}{2}}{\mid {\bf
V}\mid}^{-\frac{d-p-1}{2}}\Gamma_{p}(\frac{d-p-1}{2})},~~{\bf B}>0,
\end{equation}
where $d>2p$. This will be denoted by $IW_{p}(d,{\bf V})$.
\end{dfn}
\subsection{\bf Basic properties of elliptically contoured
distributions}~\\ We proceed here to a brief summary of the basic
properties of ECD. The following definitions and results are taken
from \cite{varga}.
\begin{dfn}
Let ${\bf X}$ be a $p\times n$ random matrix. Then, ${\bf X}$ is
said to have a matrix variate elliptically contoured distribution if
its characteristic function has the form
$$
{\phi}_{\bf X}({\bf T})~=~etr(i{\bf T}^{\prime}{\bf
M}){\psi}[tr({\bf T}^{\prime}{\bf \Sigma}{\bf T}{\bf \Phi})]
~with~{\bf T}\in{\cal R}^{p\times n},~{\bf M}\in{\cal R}^{p\times
n},~{\bf \Sigma}\in{\cal R}^{p\times p},
$$
\begin{equation}\label{charact}
{\bf \Phi}\in{\cal
R}^{n\times n},~{\bf \Sigma}~\geq 0,~~{\bf \Phi}~\geq
 0~and~{\psi}:~\lbrack~0,~\infty \lbrack~\rightarrow  R.
\end{equation}
The matrices ${\bf M}$,~${\bf \Sigma}$~and~${\bf \Phi}$~are the
parameters of the distribution.
\end{dfn}
This distribution is denoted by $ {\bf X}~\sim ~E_{p,n}({\bf M},{\bf
\Sigma}\otimes{\bf \Phi},{\psi}). $ The function ${\psi}$ is called
the {\bf characteristic generator} (cg). As a special case, when
${\psi}(z)=\exp(-\frac{z}{2})$, then ${\bf X}$ has a matrix variate
normal distribution. If $n=1$, then ${\bf x} \sim E_{p}({\bf m},
{\bf \Sigma},\psi)$ is said to have a vector variate elliptical
distribution. The relationship, in terms of the distributions, of the
matrix and the vector is illustrated as follows:
\begin{equation}\label{Aw1}
{\bf X}~\sim ~E_{p,n}({\bf M},{\bf \Sigma}\otimes{\bf \Phi},{\psi})~
if~and~only~if~{\bf x}=vec({\bf X}^{\prime})~\sim ~E_{pn}(vec({\bf
M}^{\prime}),{\bf \Sigma}\otimes{\bf \Phi},{\psi}).
\end{equation}
Here, $vec({\bf A})$ is defined by:
$$
vec({\bf A}) = \left(
\begin{array}{c}
{\bf a}_{1}\\
{\bf a}_{2}\\
\vdots\\
{\bf a}_{n}
\end{array}
\right),
$$
where ${\bf a}_{1}$,\ldots,${\bf a}_{n}$ denote the columns of the
$p\times n$ matrix ${\bf A}$. Following Schoenberg
\cite{berg}(1938), Fang and Anderson (1982)\cite{TR54} derived a
stochastic representation of matrix variate elliptical distributions
as follows.
\begin{thm}
Let ${\bf X}$ be a $p\times n$ random matrix. Let ${\bf M}$ be
$p\times n$, ${\bf \Sigma}$ be $p\times p$ and ${\bf \Phi}$ be
$n\times n$ constant matrices, with ${\bf \Sigma}\geq {\bf O}$,
${\bf \Phi}\geq {\bf O}$, $rank({\bf \Sigma})=p_{1}$, $rank({\bf
\Phi})=n_{1}$. Then
\begin{equation*}
{\bf X}~\sim ~E_{p,n}({\bf M},{\bf \Sigma}\otimes {\bf \Phi},{\psi})
\end{equation*}
if and only if
\begin{equation}\label{Aw2}
{\bf X}~\stackrel{d}{=}~{\bf M}+R{\bf A}{\bf U}{\bf B}^{\prime},
\end{equation}
where ${\bf U}$ is $p_{1}\times n_{1}$ and $vec({\bf U}^{\prime})$
is uniformly distributed on the unit sphere $S_{p_{1}n_{1}}$, $r$ is
a nonnegative random variable, $R$ and ${\bf U}$ are independent,
${\bf \Sigma}={\bf A}{\bf A}^{\prime}$ and ${\bf \Phi}={\bf B}{\bf
B}^{\prime}$ are rank factorizations of ${\bf \Sigma}$ and ${\bf
\Phi}$. Moreover
\begin{equation}\label{Aw3}
{\psi}(y)~=~\int_{0}^{\infty}{\Omega}_{p_{1}n_{1}}(R^{2}y)dF(R),~~~~y\geq
0,
\end{equation}
where ${\Omega}_{p_{1}n_{1}}({\bf t}^{\prime}{\bf t})$, ${\bf
t}~\in~{ R}^{p_{1}n_{1}}$ denotes the characteristic generator of
$vec({\bf U}^{\prime})$, and $F(r)$ denotes the distribution
function of $R$.
\end{thm}
The expression ${\bf M}+r{\bf A}{\bf U}{\bf B}^{\prime}$ is called
the stochastic representation of ${\bf X}$. In relation (\ref{Aw2}),
the notation $\stackrel{d}{=}$ stands for *equality in
distribution*. The random matrix ${\bf X}$ does not, in general,
possess a density function. But if it does, it will have the following
form (see \cite[p.26]{varga})
\begin{equation}\label{Aw4}
f({\bf X})~=~{\mid {\bf \Sigma}\mid}^{-\frac{n}{2}}{\mid {\bf
\Phi}\mid }^{-\frac{p}{2}}h(tr((\bf X-\bf M)^{\prime}{\bf
\Sigma}^{-1}(\bf X-\bf M){\bf \Phi}^{-1})),
\end{equation}
where ${\bf \Sigma}>0$ and ${\bf \Phi}>0$. The function $h$ is
called the {\bf density generator} of the distribution. In the
particular case where ${\bf \Phi}={\bf I}_{n}$ and ${\bf
M}={\bm}{\bf e}_{n}^{\prime}$, with ${\bf
e}_{n}^{\prime}=(1,1,\ldots,1)$, the pdf (\ref{Aw4}) simplifies to
the form (\ref{model1}). That is, ${\bf X}\sim E_{p,n}({\mu}{\bf
e}_{n}^{\prime},{\bf \Sigma}\otimes{\bf I}_{n},h)$, where  ${\bf
X}=({\bf x}_{1},{\bf x}_{2},\ldots,{\bf x}_{n})$ and $h(\cdot)$ is
defined through (\ref{model1}). Similarly to the family of normal
distributions, the ECD enjoy the linearity property.
\subsubsection{\bf Linearity of ECD}~\\
Let ${\bf X}~\sim ~E_{p,n}({\bf M},{\bf \Sigma}\otimes{\bf
\Phi},{\psi})$. Assume ${\bf C}$: $q\times m$, ${\bf D}$: $q\times
p$ and ${\bf B}$: $n\times m$ are constant matrices. Then
\begin{equation}\label{linear}
{\bf D}{\bf X}{\bf B}+{\bf C}~\sim ~E_{q,m}({\bf D}{\bf M}{\bf
B}+{\bf C},~({\bf D}{\bf \Sigma}{\bf D}^{\prime})\otimes({\bf
B}^{\prime}{\bf \Phi}{\bf B}),~\psi)
\end{equation}
An illustration of the linearity property is provided here:\\
Let ${\bf X}\sim E_{p,n}({\mu}{\bf e}_{n}^{\prime},{\bf
\Sigma}\otimes{\bf I}_{n},h)$ where  ${\bf X}=({\bf x}_{1},{\bf
x}_{2},\ldots,{\bf x}_{n})$. let ${\bf H}$ be the centering matrix
defined by ${\bf H}=I_{n}-\frac{1}{n}{\bf e}_{n}{\bf
e}_{n}^{\prime}$, then the matrix ${\bf H}$ is symmetric and
idempotent (that is ${\bf H}^{2}={\bf H}$). Then, the sample SSP
matrix ${\bf A}$ may be expressed as a function of ${\bf X}$:
\begin{eqnarray*}
{\bf A} & = & ({\bf X}-\overline{\bf x}{\bf e}_{n}^{\prime})({\bf X}-\overline{\bf x}{\bf e}_{n}^{\prime})^{\prime}\\
        & = & \left({\bf X}-\frac{1}{n}{\bf X}{\bf e}_{n}{\bf e}_{n}^{\prime}\right)\left({\bf X}-\frac{1}{n}{\bf X}{\bf e}_{n}{\bf e}_{n}^{\prime}\right)^{\prime} \\
        & = & {\bf X}\left({\bf I}_{n}-\frac{1}{n}{\bf e}_{n}{\bf e}_{n}^{\prime}\right)\left({\bf I}_{n}-\frac{1}{n}{\bf e}_{n}{\bf e}_{n}^{\prime}\right)^{\prime}{\bf X}^{\prime},
\end{eqnarray*}
or equivalently,
\begin{equation}\label{soul3}
{\bf A} ={\bf X}{\bf H}{\bf X}^{\prime}=({\bf X}{\bf H})({\bf X}{\bf
H})^{\prime}.
\end{equation}
Now using (\ref{linear}), we have ${\bf X}{\bf H} \sim E_{p,n}({\bf
0},{\bf \Sigma}\otimes{\bf H},h)$
\section{Some Expected Values}~\\
Estimating the precision matrix ${\bf \Sigma}^{-1}$ of the random
Kotz-Wishart matrix ${\bf A}$ requires the derivations of some
particular expected values. The following result, due to Gupta and
Varga \cite{varga}, enables us to evaluate the needed expected values of
${\bf A}$, without knowing its corresponding pdf.
\begin{thm}
Let ${\bf X}\sim E_{p,n}({\bf O},{\bf \Sigma}\otimes {\bf
\Phi},\psi)$. Let $l=rank({\bf \Sigma})$, $m=rank({\bf \Phi})$, and
$R{\bf D}{\bf U}{\bf B}^{\prime}$ be the stochastic representation
of ${\bf X}$. Assume that ${\bf Y}\sim N_{p,n}({\bf O},{\bf
\Sigma}\otimes {\bf \Phi})$. Let $K({\bf Z})$ be a function defined
on ${\cal R}^{p\times n}$ such that if ${\bf Z}\in {\cal R}^{p\times
n}$ and $a\geq 0$ then $K(a{\bf Z})=a^{k}K({\bf Z})$ where $k>-lm$.
Assume $E(K({\bf X}))$ and $E(K({\bf Y}))$ exist. Then,
\begin{equation}\label{gupta}
E(K({\bf
X}))~=~\frac{E(R^{k})\Gamma\left(\frac{lm}{2}\right)}{2^{\frac{k}{2}}\Gamma\left(\frac{lm+k}{2}\right)}E(K({\bf
Y})).
\end{equation}
\end{thm}
The proof of this theorem can be found in Gupta and Varga,
\cite[p.100]{varga}. It is important to observe that the mixing
random variables $r$ and $R$ that appear in both stochastic
representations (vector and matrix cases) are similar but not
exactly identical. Indeed, by virtue of the vec operator
(\ref{Aw1}), the moments of $R$ are obtained from those of $r$
(\ref{rkot}) by replacing $p$ by $np$. The first two moments of
${\bf A}$ are derived as follows
\begin{lemma}
Let ${\bf A}$ be the Kotz-Wishart random matrix. Then we have
\begin{equation}\label{es1}
E({\bf
A})~=~\frac{(n-1){\theta}^{-\frac{1}{s}}\Gamma\left(\frac{2q+np}{2s}\right)}{np\Gamma\left(\frac{2q+np-2}{2s}\right)}{\bf
\Sigma}=c_{1}{\bf \Sigma}.
\end{equation}
\begin{equation}\label{es2}
E\left({\bf
A}^{2}\right)~=~\frac{{\theta}^{-\frac{2}{s}}\Gamma\left(\frac{2q+np+2}{2s}\right)}{np(np+2)\Gamma\left(\frac{2q+np-2}{2s}\right)}E\left({\bf
W}^{2}\right),
\end{equation}
where ${\bf W}\stackrel{d}{=}{\bf W}_{p}({\bf \Sigma},n-1)$ and
$E\left({\bf W}^{2}\right)$ is given by (see \cite[p.4]{adhi})
\begin{equation}\label{afrs1}
E\left({\bf W}^{2}\right)~=~(n-1)^{2}{\bf
\Sigma}^{2}+(n-1)\left[{\bf \Sigma}tr({\bf \Sigma})+{\bf
\Sigma}^{2}\right].
\end{equation}
\end{lemma}
{\bf Proof:} Let us define the function $K({\bf X})~=~{\bf
X}\left({\bf I}_{n}-\frac{1}{n}{\bf e}_{n}{\bf
e}_{n}^{\prime}\right){\bf X}^{\prime}~=~{\bf X}{\bf H}{\bf
X}^{\prime}$. Hence, $E(K({\bf X}))=E({\bf A})$, and $E(K({\bf
Y}))=E({\bf W}_{p}({\bf \Sigma},n-1))=(n-1){\bf \Sigma}$. Using
equation (\ref{gupta}), where $l=p, m=n$, $k=2$, and equation
(\ref{rkot}), (with $p$
replaced by np),we obtain the desired result (\ref{es1}).\\
 To prove the equation
(\ref{es2}), we consider the function $K(\cdot)$ defined on
${R}^{p\times n}$ by
\begin{equation*}
K({\bf X})~=~\left[{\bf X}\left({\bf I}_{n}-\frac{1}{n}{\bf
e}_{n}{\bf e}_{n}^{\prime}\right){\bf X}^{\prime}\right]^{2}.
\end{equation*}
So, we have $E(K({\bf X}))=E\left({\bf A}^{2}\right)$, and for a given $a\geq 0$ $K(a{\bf X})=a^{4}K({\bf X})$.\\
Now, using equation \ref{gupta} (with $k=4$) and the fact that
$K({\bf Y})\stackrel{d}{=}{\bf W}^{2}$, we can write
\begin{equation}\label{om1}
E\left({\bf
A}^{2}\right)~=\frac{E(r^4)\Gamma\left(\frac{np}{2}\right)}{2^{2}\Gamma\left(\frac{np+4}{2}\right)}{\cdot}E\left({\bf
W}^{2}\right).
\end{equation}
From equation \ref{gupta} (with $t=2$ and $p$ replaced by $np$), we
can also write
\begin{equation}\label{om2}
E\left(R^{4}\right)~=~\frac{{\theta}^{-\frac{2}{s}}\Gamma\left(\frac{2q+np+2}{2s}\right)}{\Gamma\left(\frac{2q+np-2}{2s}\right)}.
\end{equation}
Consequently, we have
\begin{eqnarray*}
E\left({\bf A}^{2}\right) & = & \frac{{\theta}^{-\frac{2}{s}}\Gamma\left(\frac{2q+np+2}{2s}\right)\Gamma\left(\frac{np}{2}\right)}{4\Gamma\left(\frac{2q+np-2}{2s}\right)\Gamma\left(\frac{np+4}{2}\right)}E\left({\bf W}^{2}\right)\\
                          & = & \frac{{\theta}^{-\frac{2}{s}}\Gamma\left(\frac{2q+np+2}{2s}\right)}{np(np+2)\Gamma\left(\frac{2q+np-2}{2s}\right)}{\cdot}E\left({\bf W}^{2}\right),
\end{eqnarray*}
which completes the proof of the lemma.
\begin{lemma}
Let ${\bf A}$ be the Kotz-Wishart random matrix. Then for any $t>0$
the t th moment of the generalized variance $\mid {\bf A}\mid $ is
given by
\begin{equation}\label{aw3}
E\left(\mid{\bf A}\mid
^{t}\right)~=~\frac{{\theta}^{-\frac{tp}{s}}\Gamma\left(\frac{2q+np+2tp-2}{2s}\right)\Gamma\left(\frac{np}{2}\right){\Gamma}_{p}\left(\frac{n-1}{2}+t\right)}{\Gamma\left(\frac{2q+np-2}{2s}\right)\Gamma\left(\frac{np+2tp}{2}\right){\Gamma}_{p}\left(\frac{n-1}{2}\right)}\mid
{\bf \Sigma}\mid ^{t}.
\end{equation}
\end{lemma}
{\bf Proof:} To prove (\ref{aw3}), we consider the function
$K(\cdot)$ defined on ${{R}}^{p\times n}$ by
\begin{equation*}
K({\bf X})~=~\mid {\bf X}({\bf I}_{n}-\frac{1}{n}{\bf e}_{n}{\bf
e}_{n}^{\prime}){\bf X}^{\prime}\mid ^{t},~~t>0~.
\end{equation*}
Hence, $E(K({\bf X}))=E\left(\mid {\bf A}\mid ^{t}\right)$ and for a
given $a\geq 0$ we have
\begin{equation*}
K(a{\bf X})~=~a^{2tp}K({\bf X}).
\end{equation*}
By using equation \ref{gupta}, (with $k=2tp$), we can write
\begin{equation}\label{gup}
E\left(\mid {\bf A}\mid
^{t}\right)~=~\frac{E\left(r^{2tp}\right)\Gamma\left(\frac{np}{2}\right)}{2^{tp}\Gamma\left(\frac{np+2tp}{2}\right)}{\cdot}E(K({\bf
Y})).
\end{equation}
Here $E(K({\bf Y}))=E\left(\mid {\bf W}\mid ^{t}\right)$ where ${\bf
W}\stackrel{d}{=}{\bf W}_{p}({\bf \Sigma},n-1)$, (that is a Wishart
distribution). So, $E\left(\mid {\bf W}\mid ^{t}\right)$ is given by
the following well-known result (see \cite{asarr} or \cite[p.101]{muir})
\begin{equation}\label{mred}
E\left(\mid {\bf W}\mid
^{t}\right)~=~\frac{2^{tp}{\Gamma}_{p}\left(\frac{n-1}{2}+t\right)}{{\Gamma}_{p}\left(\frac{n-1}{2}\right)}{\cdot}\mid
{\bf \Sigma}\mid ^{t},
\end{equation}
and $E(r^{2tp})$ is given by (see equation \ref{rkot})
\begin{equation}
E\left(R^{2tp}\right)~=~\frac{{\theta}^{-\frac{tp}{s}}\Gamma\left(\frac{2q+np+2tp-2}{2s}\right)}{\Gamma\left(\frac{2q+np-2}{2s}\right)}.
\end{equation}
Now, substituting (\ref{rkot}) and (\ref{mred}) in equation (\ref{gup}), we obtain the desired result (\ref{aw3}).
For the calculations of $E({\bf A}^{-1}), E({\bf A}^{-2})$ and
$E(\mid{\bf A}\mid^{-t})$, interested readers are referred to
\cite{sarr}.
\section{Estimation of the Precision Matrix}~\\
The problem of estimating the precision matrix ${\bf \Sigma}^{-1}$
of a multivariate normal model has been widely investigated by
several authors. A review paper on this particular subject was
presented by Kubokawa \cite{kub}. Specially in an empirical Bayes
estimation context, Efron and Morris \cite{efron} had proved that
the unbiased estimator of the precision matrix was the best constant
multiple of ${\bf A}^{-1}$. Hence, the purpose of this section
consists in proving that the above results remain robust under the
multivariate Kotz model. To achieve our goal here, we will make use
of the same Efron-Morris's loss function, specified by:
\begin{equation}\label{EfronM}
l({\bf \Sigma}^{-1},{\bf \Delta}) = \frac{tr\left[({\bf \Delta}-{\bf
\Sigma}^{-1})^{2}{\bf A}\right]}{\nu tr({\bf \Sigma}^{-1})},
\end{equation}
where ${\bf \Delta}$ denotes any estimator of ${\bf \Sigma}^{-1}$
and ${\bf A}$ is $W_{^p}(\nu,\Sigma)$, with $\nu >p+1$. As usual,
the corresponding risk function is defined by
\begin{equation}\label{riskf}
R({\bf \Sigma}^{-1},{\bf \Delta})=E\left[l({\bf \Sigma}^{-1},{\bf
\Delta})\right],
\end{equation}
where the expectation is taken with respect to the distribution of
the random matrix ${\bf A}$. In our present case, ${\bf A}
\stackrel{d}{=} KW_{p}(n-1,\Sigma)$, and we will assume $n>p+2$.
Therefore, some adjustments on the Efron-Morris's loss function are
made accordingly; say $\nu=n-1$. Note that the problem of estimating
${\bf \Sigma}^{-1}$ was considered in \cite{sarr}, under a quadratic
loss function. Further, an unbiased estimator ${\bf
\Delta}_{0}=c_{0}{\bf A}^{-1}$ for the precision matrix was
obtained, where
\begin{equation}\label{unbias}
c_{0}=\frac{(n-p-2)\Gamma
\left(\frac{2q+np-2}{2s}\right)}{{\theta}^{1/s}(np-2)\Gamma
\left(\frac{2q+np-4}{2s}\right)}.
\end{equation}
Recall that an estimator ${\bf \Delta}_{1}$ is said to be better
than another estimator ${\bf \Delta}_{2}$ if $R({\bf
\Sigma}^{-1},{\bf \Delta}_{1})\leq R({\bf \Sigma}^{-1},{\bf
\Delta}_{2}),~~\forall$ ${\bf \Sigma}>0$. The dominance result involving
the unbiased estimator ${\bf \Delta}_{0}$ is derived in the
following theorem.
\begin{thm}~\\
The best (smallest risk) estimator of the precision matrix ${\bf
\Sigma}^{-1}$ of the Kotz-Wishart distribution, having the form $\alpha
{\bf A}^{-1}$, is the unbiased estimator ${\bf \Delta}_{0}=c_{0}{\bf
A}^{-1}$, where $c_{0}$ is given by (\ref{unbias}). And its
corresponding risk is given by
\begin{equation}\label{risk}
R({\bf \Sigma}^{-1},{\bf \Delta}_{0})=\frac{c_{1}-c_{0}}{\nu},
\end{equation}
where $c_{1}$ is specified in (\ref{es1}) and $\nu = n-1$.
\end{thm}
{\bf Proof:}~\\
Let ${\bf \Delta}=\alpha {\bf A}^{-1}$, $\alpha >0$ be an estimator
of ${\bf \Sigma}^{-1}$. Since $E(tr(\cdot))=tr(E(\cdot))$, and both
functions $E(\cdot)$, $tr(\cdot)$ are linear, we have:
\begin{eqnarray*}
{\nu}tr({\bf \Sigma}^{-1})R({\bf \Sigma}^{-1},{\bf \Delta}) & = &
E\left[{\alpha}^{2}tr({\bf A}^{-1})-2{\alpha}tr({\bf
\Sigma}^{-1})+tr({\bf \Sigma}^{-2}{\bf A})\right] \\        & = &
{\alpha}^{2}\frac{1}{c_{0}}tr({\bf \Sigma}^{-1})-2{\alpha}tr({\bf
\Sigma}^{-1})+c_{1}tr({\bf \Sigma}^{-1})\\                  & = &
\left(\frac{{\alpha}^{2}}{c_{0}}-2 \alpha+c_{1}\right)tr({\bf
\Sigma}^{-1}).
\end{eqnarray*}
Thus, $R({\bf \Sigma}^{-1},{\bf \Delta})=\dfrac{g(\alpha)}{\nu}$,
where $g(\alpha)=\frac{{\alpha}^{2}}{c_{0}}-2 \alpha+c_{1}$ is a
quadratic function in $\alpha$, which attains its minimum at
${\alpha}_{0}$ solution of the equation $g^{\prime}(\alpha)=0$, here
$g^{\prime}(\cdot)$ denotes the derivative. Hence, the unique
solution is ${\alpha}_{0}=c_{0}$, which proves that ${\bf
\Delta}_{0}$ is the best estimator of ${\bf \Sigma}^{-1}$, among the
estimators of the form ${\alpha}{\bf A}^{-1}$. Its corresponding
risk is:
\begin{eqnarray*}
R({\bf \Sigma}^{-1},{\bf \Delta}_{0}) & = &\frac{g(c_{0})}{\nu}\\
                                      & = &\frac{c_{1}-c_{0}}{\nu},
\end{eqnarray*}
which completes the proof of the theorem.
\begin{rmq}~\\As a special case, when $q=2\theta=s=1$, then the Kotz
model reduces to the usual classical normal model. In this case, the
constants $c_{0}$ and $c_{1}$ reduce to $n-p-2$ and $n-1$
respectively. Consequently, the above risk simplifies to
\begin{equation*}
R({\bf \Sigma}^{-1},{\bf \Delta}_{0})=\frac{p+1}{\nu},
\end{equation*}
which is the risk obtained by Efron and Morris\cite{efron}. Outside
the class of constant multiples of ${\bf A}^{-1}$, the quoted authors
proposed better estimators than the unbiased one. Instead of
following their footsteps along that direction,
we just recall that similar improved estimators have been considered
in Sarr and Gupta \cite{sarr}, under a quadratic loss function.
Recently, increasing interest in the problem of estimating the
precision matrices is observed in high dimensional contexts; with
applications in discriminant analysis \cite{sriv} as well as in
finance \cite{bai}.
\end{rmq}
\section{Derivation of Density Formulas}~\\
In this section, we provide the explicit expressions for the
probability density functions for the Kotz-Wishart random matrix and
for its inverted version.\\ For any uncorrelated (or dependent)
elliptical model of the form (\ref{model1}), Anderson and Fang (see
\cite[p.207]{gth} or \cite[p.239]{varga}) derived a general formula
for the pdf of its corresponding SSP matrix ${\bf A}$:
\begin{equation}\label{genwish}
f({\bf A})~=~\frac{2{\pi}^{\frac{pn}{2}}\mid {\bf A}\mid
^{\frac{n-p}{2}-1}}{\Gamma\left(\frac{p}{2}\right)\Gamma_{p}\left(\frac{n-1}{2}\right)\mid
{\bf \Sigma}\mid
^{\frac{n-1}{2}}}\int_{0}^{\infty}x^{p-1}h(x^{2}+tr({\bf
\Sigma}^{-1}{\bf A}))dx.
\end{equation}
Here, we derive an explicit form of the pdf of the Kotz-Wishart
random matrix ${\bf A}$, when $s=1$. In other words, we have to
compute the above integral when
\begin{equation}\label{ous1}
h(y)~=~c_{pn}y^{q-1}\exp(-\theta y),~~~y\geq
0~and~~c_{pn}~=~\frac{{\theta}^{\frac{2q+pn-2}{2}}\Gamma\left(\frac{pn}{2}\right)}{{\pi}^{\frac{pn}{2}}\Gamma\left(\frac{2q+pn-2}{2}\right)}~.
\end{equation}
The result is given in the following lemma
\begin{lemma}
The pdf of the Kotz-Wishart random matrix ${\bf A}$ is
\begin{eqnarray}\label{genkotz}
f({\bf A}) & = & K\mid {\bf A}\mid ^{\frac{n-p}{2}-1}\left(tr({\bf \Sigma}^{-1}{\bf
A})\right)^{\frac{2q+p-4}{4}}\exp\left(-\frac{\theta tr({\bf
\Sigma}^{-1}{\bf A})}{2}\right)\nonumber \\
           &   & \times W_{\alpha,\beta}(\theta tr({\bf \Sigma}^{-1}{\bf A})),
\end{eqnarray}
where
\begin{equation}\label{coa}
K~=~\frac{\Gamma\left(\frac{pn}{2}\right)\mid {\bf \Sigma}\mid
^{-\frac{n-1}{2}}}{{\Gamma}_{p}\left(\frac{n-1}{2}\right)\Gamma\left(\frac{2q+pn-2}{2}\right)}{\theta}^{\frac{2q+2pn-p-4}{4}}~,
\end{equation}
and $W_{\alpha,\beta}(\cdot)$ denotes the Whittaker function (see
the Appendix), with
\begin{equation}\label{indwit}
\alpha~=~\frac{2q-p}{4}~,~~~~\beta~=~\frac{2q+p-2}{4}~.
\end{equation}
\end{lemma}
{\bf Proof:} Let us calculate the integral
$$
{\bf I} = \int_{0}^{\infty}x^{p-1}h\left(x^2+tr({\bf \Sigma}^{-1}{\bf
A})\right)dx,
$$
where $h(\cdot)$ is given by (\ref{ous1}). Setting $b=tr({\bf
\Sigma}^{-1}{\bf A})$ leads to
\begin{eqnarray*}
{\bf I} & = & c_{pn}\int_{0}^{\infty}x^{p-1}(x^2+b)^{q-1}\exp\left(-\theta (x^2+b)\right)dx \\
        & = & c_{pn}\exp(-\theta b)\int_{0}^{\infty}x^{p-1}(x^2+b)^{q-1}\exp(-\theta{x^2})dx.
\end{eqnarray*}
Let us make the change of variable $z=x^2$, then the jacobian of the
transformation is $\frac{1}{2}z^{-\frac{1}{2}}$. Hence we obtain
\begin{eqnarray*}
\int_{0}^{\infty}x^{p-1}(x^2+b)^{q-1}\exp(-\theta {x^2})dx & = & \int_{0}^{\infty}z^{\frac{p-1}{2}}(z+b)^{q-1}\exp(-\theta z)\frac{1}{2}z^{-\frac{1}{2}}dz \\
                                                           & = & \frac{1}{2}\int_{0}^{\infty}z^{\frac{p}{2}-1}(z+b)^{q-1}\exp(-\theta z)dz .
\end{eqnarray*}
Now using equation (\ref{mellin}), taken from the appendix, we
can write
\begin{equation*}
\int_{0}^{\infty}z^{\frac{p}{2}-1}(z+b)^{q-1}\exp(-\theta
z)dz~=~b^{\frac{2q+p-4}{4}}{\theta}^{-\frac{2q+p}{4}}\Gamma\left(\frac{p}{2}\right)\exp\left(\frac{\theta
b}{2}\right)W_{\alpha,\beta}(\theta b)~,
\end{equation*}
where
$$
\alpha~=~\frac{2q-p}{4},~~~\beta~=~\frac{2q+p-2}{4}.
$$
Thus
\begin{equation*}
\int_{0}^{\infty}x^{p-1}h(x^2+b)dx~=~\frac{c_{pn}}{2}{\theta}^{-\frac{2q+p}{4}}\Gamma\left(\frac{p}{2}\right)b^{\frac{2q+p-4}{4}}\exp\left(-\frac{\theta
b}{2}\right)W_{\alpha ,\beta}(\theta b)~.
\end{equation*}
Define $K_{1}$ as follows
\begin{equation*}
K_{1}~=~\frac{c_{pn}}{2}\Gamma\left(\frac{p}{2}\right){\theta}^{-\frac{2q+p}{4}}~=~\frac{\Gamma\left(\frac{pn}{2}\right)\Gamma\left(\frac{p}{2}\right)}{2{\pi}^{\frac{pn}{2}}\Gamma\left(\frac{2q+np-2}{2}\right)}{\theta}^{\frac{2q+2pn-p-4}{4}},
\end{equation*}
we have
\begin{equation*}
\frac{2{\pi}^{\frac{pn}{2}}\mid {\bf \Sigma}\mid
^{-\frac{n-1}{2}}}{\Gamma\left(\frac{p}{2}\right){\Gamma}_{p}\left(\frac{n-1}{2}\right)}K_{1}~=~\frac{\Gamma\left(\frac{pn}{2}\right)\mid
{\bf \Sigma}\mid
^{-\frac{n-1}{2}}}{{\Gamma}_{p}\left(\frac{n-1}{2}\right)\Gamma\left(\frac{2q+pn-2}{2}\right)}{\theta}^{\frac{2q+2pn-p-4}{4}}~=K,
\end{equation*}
where $K$ is given by (\ref{coa}). So, we finally get
\begin{equation*}
\frac{2{\pi}^{\frac{pn}{2}}\mid {\bf A}\mid
^{\frac{n-p}{2}-1}}{\Gamma\left(\frac{p}{2}\right)\Gamma_{p}\left(\frac{n-1}{2}\right)\mid
{\bf \Sigma}\mid
^{\frac{n-1}{2}}}\int_{0}^{\infty}x^{p-1}h(x^{2}+tr({\bf
\Sigma}^{-1}{\bf A}))dx~=~
\end{equation*}
\begin{equation*}
K \mid {\bf A}\mid ^{\frac{n-p}{2}-1}\left(tr({\bf \Sigma}^{-1}{\bf
A})\right)^{\frac{2q+p-4}{4}}\exp\left(-\frac{\theta tr({\bf
\Sigma}^{-1}{\bf A})}{2}\right)W_{\alpha,\beta}(\theta tr({\bf
\Sigma}^{-1}{\bf A})),
\end{equation*}
which completes the proof of the lemma.\\
\begin{rmq}\label{mark}
As a special case, if $q=2\theta =1$, then $\alpha$ and $\beta$
reduce to $\frac{2-p}{4}$ and $\frac{p}{4}$ respectively. Now, using
the identity (\ref{simpW}) from the Appendix enables us to simplify
the pdf (\ref{genkotz}):
\begin{equation}
f({\bf A})~=~\frac{\mid {\bf A}\mid
^{\frac{n-p}{2}-1}}{2^{\frac{p(n-1)}{2}}{\Gamma}_{p}\left(\frac{n-1}{2}\right)\mid
{\bf \Sigma}\mid ^{\frac{n-1}{2}}}\exp\left(-\frac{1}{2}tr({\bf
\Sigma}^{-1}{\bf A})\right),
\end{equation}
which is identical to (\ref{Wishart}).That is, ${\bf
A}~\sim~W_{p}(n-1,\Sigma)$. The identity (\ref{simpW}) will be used repeatedly in the present work.
\end{rmq}
The expression (\ref{genkotz}) may also be rewritten as
\begin{eqnarray}\label{simpdf}
f({\bf A}) & = & C_{1}(q,\theta)\mid{\bf \Sigma}\mid
^{-\frac{\nu}{2}}\mid{\bf A}\mid ^{\frac{\nu-p-1}{2}}\left(\theta
tr({\bf \Sigma}^{-1}{\bf A})\right)^{\frac{2q+p-4}{4}}\exp\left(-\frac{\theta tr({\bf
\Sigma}^{-1}{\bf A})}{2}\right)\nonumber \\
           &   & \times W_{\alpha,\beta}(\theta tr({\bf \Sigma}^{-1}{\bf A})),
\end{eqnarray}
where $\nu=n-1$, $\alpha=\frac{2q-p}{4}, \beta=\frac{2q+p-2}{4}$ and
\begin{equation*}
C_{1}(q,\theta)=\frac{\Gamma\left(\frac{p(\nu+1)}{2}\right){\theta}^{\frac{p\nu}{2}}}{\Gamma\left(\frac{2q+p(\nu+1)-2}{2}\right){\Gamma}_{p}(\frac{\nu}{2})}.
\end{equation*}
\begin{dfn}~\\
Under the restriction $s=1$, a $p\times p$ symmetric random matrix
${\bf A}>0$ is said to have a Kotz-Wishart distribution, with $\nu$
degrees of freedom and parameter matrix ${\bf \Sigma}>0$ if its pdf
is given by (\ref{simpdf}), with $\nu \geq p$.
\end{dfn}
The pdf of the random matrix ${\bf B}={\bf A}^{-1}$ is derived in
the next lemma
\begin{lemma}~\\
Let ${\bf A}$ be the Kotz-Wishart matrix whose pdf is
(\ref{simpdf}). Let ${\bf B}={\bf A}^{-1}$. Then the pdf of ${\bf
B}>0$ is
\begin{eqnarray}\label{invKotz}
g({\bf B})& = & C_{1}\mid{\bf \Sigma}\mid ^{-\frac{\nu}{2}}\mid{\bf B}\mid ^{-\frac{\nu+p+1}{2}}\left(\theta
tr({\bf \Sigma}^{-1}{\bf B}^{-1})\right)^{\frac{2q+p-4}{4}}\exp\left(-\frac{\theta tr({\bf
\Sigma}^{-1}{\bf B}^{-1})}{2}\right)\nonumber \\
          &   & \times W_{\alpha,\beta}(\theta tr({\bf \Sigma}^{-1}{\bf B}^{-1})),
\end{eqnarray}
where $\nu, \alpha, \beta$ and $C_{1}=C_{1}(q,\theta)$ are given in
(\ref{simpdf}).
\end{lemma}
{\bf Proof:}~\\
Since the jacobian of the transformation ${\bf B}={\bf A}^{-1}$ is
$\mid {\bf B} \mid ^{-(p+1)}$, (see \cite[p.59]{muir}), the desired
result follows directly.
\begin{dfn}{\bf Inverted Kotz-Wishart}~\\
A $p\times p$ random matrix ${\bf B}>0$ is said to have an inverted
Kotz-Wishart distribution, with d degrees of freedom and parameter
matrix ${\bf V}>0$ if its pdf is
\begin{eqnarray}\label{invpdf}
g({\bf B})& = & C_{1}(q,\theta)\mid{\bf V}\mid ^{\frac{d-p-1}{2}}\mid{\bf
B}\mid ^{-\frac{d}{2}}\left(\theta tr({\bf V}{\bf
B}^{-1})\right)^{\frac{2q+p-4}{4}}\exp\left(-\frac{\theta tr({\bf
 V}{\bf B}^{-1})}{2}\right)\nonumber \\
          &   & \times W_{\alpha,\beta}(\theta tr({\bf V}{\bf B}^{-1})),
 \end{eqnarray}
 where $d>2p$, and $\alpha, \beta$ given in (\ref{simpdf}). This will be denoted by ${\bf B}$ is $IKW_{p}(d,{\bf V})$.
 \end{dfn}
 Once again, when $q=2\theta=1$, then $\alpha+\beta=\frac{1}{2}$ and by
 virtue of (\ref{simpW}), the above pdf reduces to
 \begin{equation*}
 g({\bf
 B})=\frac{2^{-\frac{p(d-p-1)}{2}}}{{\Gamma}_{p}\left(\frac{d-p-1}{2}\right)}\mid
 {\bf V} \mid ^{\frac{d-p-1}{2}}\mid {\bf B}\mid
 ^{-\frac{d}{2}}etr\left(-\frac{1}{2}{\bf V}{\bf
 B}^{-1}\right),~~~{\bf B}>0,
 \end{equation*}
 which coincides with (\ref{invWishart}), the pdf of the inverted
 Wishart distribution.\\
 A direct consequence of the above definitions is illustrated
 through the following relationship:
 $$
 {\bf A} \stackrel{d}{=} KW_{p}(\nu, \Sigma)~~\Rightarrow
{\bf A}^{-1} \stackrel{d}{=} IKW_{p}(\nu+p+1, {\Sigma}^{-1}).
$$
 Further properties of the Kotz-Wishart distributions as well as its
 inverted version will be investigated in the next section.
 \section{Some Results Involving Zonal Polynomials}~\\
 Zonal polynomials play a crucial role in multivariate
 statistical analysis. Indeed, the exact distributions of many
 statistics involve, in general, zonal polynomials, as pointed out
 in the excellent textbook by Muirhead \cite{muir}. The theory of
 zonal polynomials for real matrices was first introduced by James
 (1961) \cite{james61}. More generally, the foundations of the
 theory of zonal polynomials and its applications in multivariate
 analysis were built through a series of papers \cite{james60},
 \cite{james61}, \cite{const} and \cite{james64} by James and
 Constantine. And most of the main results involving zonal
 polynomials are collected in Muirhead's book. However, those existing
 results are all essentially based on the normality assumption.
 Generalized multivariate analysis, our current framework, involves more general distributions than the normal
 ones, such as elliptically contoured distributions (ECD). Before
 deriving our main results, let us first recall some useful
 notations and definitions. For more details about the construction
 of zonal polynomials, interested readers are referred to
 \cite{muir}, chapter.7.\\ The Pochammer symbol $(a)_{k}$ is defined by
 $(a)_{k} = a(a+1)\ldots (a+k-1),~k=1,2\ldots $ with $(a)_{0}=1$.
 Let $C_{\kappa}({\bf X})$ be the zonal polynomial of the $p\times p$
 symmetric matrix ${\bf X}$ corresponding to the ordered partition
 $\kappa = (k_{1},\ldots, k_{p}),~k_{1}+k_{2}+\ldots +k_{p} = k $,
 $k_{1}\geq k_{2}\geq \ldots \geq k_{p}\geq 0$. Then
 \begin{equation}\label{trace}
 (tr{\bf X})^{k} = \sum_{\kappa}C_{\kappa}({\bf X}),
 \end{equation}
 where $\sum_{\kappa}$ denotes summation over all the ordered
 partitions $\kappa$ of k. The generalized hypergeometric functions, which involve zonal
 polynomials, are defined as follows (see [p.258-259]\cite{muir}):
 \begin{dfn}~The hypergeometric functions of matrix argument are
 given by
 \begin{equation}\label{hypgeo}
 _{m}F_{q}(a_{1},\ldots, a_{m};b_{1},\ldots,b_{q};{\bf X}) =
 \sum_{k=0}^{\infty}\sum_{\kappa}\frac{(a_{1})_{\kappa}\cdots
 (a_{m})_{\kappa}}{(b_{1})_{\kappa} \cdots (b_{q})_{\kappa}}\frac{C_{\kappa}({\bf
 X})}{k!},
 \end{equation}
where $a_{1},\ldots, a_{m},b_{1},\ldots, b_{q}$ are real or complex
and the generalized hypergeometric coefficient $(a)_{\kappa}$ is
given by
\begin{equation}
(a)_{\kappa} =\prod_{i=1}^{p}(a-\frac{1}{2}(i-1))_{k_{i}},
\end{equation}
where $(r)_{k_{i}}$ denotes the Pochammer symbol. The matrix ${\bf
X}$ being a complex symmetric $ p\times p$ matrix, and it is also assumed that $m\leq q+1$.
\end{dfn}
Two special cases, which are of great interest in our work, are
\begin{equation}\label{F00}
_{0}F_{0}({\bf X})
=\sum_{k=0}^{\infty}\sum_{\kappa}\frac{C_{\kappa}({\bf
 X})}{k!}=\sum_{k=0}^{\infty}\frac{(tr{\bf X})^{k}}{k!} = etr({\bf X}).
\end{equation}
 And
\begin{equation}\label{F}
 _{1}F_{0}(a;{\bf X})=\sum_{k=0}^{\infty}\sum_{\kappa}(a)_{\kappa}\frac{C_{\kappa}({\bf
 X})}{k!}=|I_{p}-{\bf X}|^{-a},~~~~~~||{\bf X}||<1,
\end{equation}
 where $||{\bf X}|| $ denotes the maximum of the absolute values of the eigenvalues of ${\bf X}$.
\begin{rmq}\label{Fpoly}~\\
If any numerator parameter $a_{1}$ is a negative integer, say, $a_{1}=-m $, then the function (\ref{hypgeo}) reduces to a polynomial of degree $ pm $, because for $k\geq pm+1$, $(a_{1})_{\kappa}=0$. Khatri (\cite{khatri}) exploited this particular case to derive explicit forms of the distributions of the smallest eigenvalue of a Wishart matrix. Here, we will establish similar results based on the Kotz-Wishart matrix.
\end{rmq}
First let us review some useful lemmas, needed to achieve our goal.
\subsection{\bf Some preliminary results on integration}~\\
The following lemma, due to Li (1997) (\cite{rli}), is completely proved in (\cite{diaz}) (2011).
\begin{lemma}~\\
let ${\bf Z}$ be a complex symmetric $p\times p$ matrix with $Re({\bf Z})>0$, and let ${\bf U}$ be a symmetric $p\times p$ matrix. Assume
$$
\gamma =\int_{0}^{\infty}f(y)y^{pa-k-1}dy < \infty ,
$$
where $k$ is a given positive integer. Then
\begin{equation}\label{LiR}
\int_{{\bf X} >0}|{\bf X}|^{a-\frac{p+1}{2}}f(tr{\bf X}{\bf Z})C_{\kappa}({\bf
X}^{-1}{\bf U})(d{\bf X})=
\end{equation}
\begin{equation*}
\frac{(-1)^{k}\Gamma_{p}(a)}{\Gamma(ap-k)(-a+\frac{p+1}{2})_{\kappa}}{\cdot}{\gamma}{\cdot} |{\bf Z}|^{-a}C_{\kappa}({\bf
U}{\bf Z}),
\end{equation*}
for $Re(a) > k_{1}+\frac{p-1}{2}$, where $\kappa = (k_{1},\ldots, k_{p})$, and $k_{1}+k_{2}+\ldots +k_{p} = k $.
\end{lemma}
The following lemma, presented in (\cite{caro}), is needed to derive the moment generating function (mgf) of the Kotz-Wishart matrix.
\begin{lemma}~\\
let ${\bf Z}$ be a complex symmetric $p\times p$ matrix with $Re({\bf Z})>0$, and let ${\bf U}$ be a symmetric $p\times p$ matrix. Assume
$$
\delta =\int_{0}^{\infty}h(y)y^{pa+k-1}dy < \infty ,
$$
where $k$ is a given positive integer. Then
\begin{equation}\label{garcia}
\int_{{\bf X} >0}|{\bf X}|^{a-\frac{p+1}{2}}h(tr{\bf X}{\bf Z})C_{\kappa}({\bf
X}{\bf U})(d{\bf X})=\frac{(a)_{\kappa} \Gamma_{p}(a)}{\Gamma(ap+k)}{\cdot}{\delta}{\cdot} |{\bf Z}|^{-a}C_{\kappa}({\bf
U}{\bf Z}^{-1}),
\end{equation}
for $Re(a) > \frac{p-1}{2}$, where $\kappa = (k_{1},\ldots, k_{p})$, and $k_{1}+k_{2}+\ldots +k_{p} = k $.
\end{lemma}
The next result will be used to calculate the constant $\gamma$ in (\ref{LiR}).
\begin{lemma}~\\
Let ${\bf \Sigma}>0$ and ${\bf \Lambda}>0$ be two $p\times p$ matrices. Let $h(\cdot)$ be defined as follows:
\begin{eqnarray*}
h(tr{\bf \Sigma}^{-1}{\bf X}) & = & h_{1}(tr{\bf \Sigma}^{-1}{\bf X}+d)\nonumber \\
                              & = & \left[\theta(tr{\bf \Sigma}^{-1}{\bf X}+d)\right]^{\xi}\exp(-\frac{\theta}{2}(tr{\bf \Sigma}^{-1}{\bf X}+d))W_{\alpha,\beta}(\theta(tr{\bf \Sigma}^{-1}{\bf X}+d)),
\end{eqnarray*}
where ${\bf X}>0$ is a $p\times p$ matrix, $d=tr{\bf \Sigma}^{-1}{\bf \Lambda}$, $\theta >0$, $q>0$, $\alpha=\frac{2q-p}{4}$, $\beta=\frac{2q+p-2}{4}$ and $\xi=\frac{2q+p-4}{4}$. Then
\begin{equation}\label{mar1}
\int_{0}^{\infty}h(z)z^{\frac{\nu p}{2}-k-1}dz=
\end{equation}
\begin{equation*}
\Gamma\left(\frac{\nu p}{2}-k\right)\left(tr{\bf \Sigma}^{-1}{\bf \Lambda}\right)^{\frac{\nu p}{2}-k}{\bf
G}^{30}_{23}\left({\theta}tr{\bf \Sigma}^{-1}{\bf
\Lambda}\begin{tabular}{|ccc}
                           $ 0$,                & $\frac{p}{2}$      &  \\
                           $-\frac{p\nu}{2}+k$, & $q-1+\frac{p}{2}$, & 0 \\
                     \end{tabular} \right)~,
\end{equation*}
where k is a given positive integer and $G^{ts}_{qr}$ denotes the Meijer's ${\bf G}$-function
\end{lemma}
{\bf Proof:}~\\
\begin{eqnarray*}
\gamma & = & \int_{0}^{\infty}h(z)z^{\frac{\nu p}{2}-k-1}dz \\
       & = & \int_{0}^{\infty}[\theta (z+d)]^{\frac{2q+p-4}{4}}\exp(-\frac{\theta}{2}(z+d))W_{\alpha,\beta}(\theta(z+d))z^{\frac{\nu p}{2}-k-1}dz.
\end{eqnarray*}
Using the change of variable ${\theta}z=t$ gives:
\begin{eqnarray*}
\gamma & = & e^{-\frac{\theta d}{2}}\int_{0}^{\infty}(t+\theta d)^{\frac{2q+p-4}{4}}e^{-\frac{t}{2}}W_{\alpha,\beta}(t+\theta d){\theta}^{-\frac{\nu p}{2}+k}t^{\frac{\nu p}{2}-k-1}dt\\
       & = & e^{-\frac{\theta d}{2}}{\theta}^{-\frac{\nu p}{2}+k}\int_{0}^{\infty} t^{\frac{\nu p}{2}-k-1}(t+\theta d)^{\frac{2q+p-4}{4}}e^{-\frac{t}{2}}W_{\alpha,\beta}(t+\theta d)dt
\end{eqnarray*}
The last integral is evaluated using the identity (\ref{Meijer}), by substituting the appropriate parameters. Namely, $\rho=\frac{\nu p}{2}-k$, $c= \theta tr{\bf \Sigma}^{-1}{\bf \Lambda}$, $\sigma=-\frac{2q+p-4}{4}$, $\alpha=\frac{2q-p}{4}$, $\beta=\frac{2q+p-2}{4}$, $1-\alpha-\sigma =\frac{p}{2}$, $\frac{1}{2}+\beta-\sigma=q-1+\frac{p}{2}$, and $\frac{1}{2}-\beta-\sigma=0$.\\
With the help of the above materials, we are now ready to extend Khatri's results (\cite{khatri}).
\subsection{\bf Distribution of the smallest eigenvalue of KW matrix}~\\
The distribution of the smallest eigenvalue of the Kotz-Wishart matrix ${\bf A}$ is obtained by following Khatri's approach: First, we evaluate the
probability $P({\bf A} > {\bf \Lambda})$, where ${\bf \Lambda}>0$ and ${\bf A} > {\bf \Lambda}$ means ${\bf A}-{\bf \Lambda}>0$ . Secondly, we use the obtained probability to derive the cdf of the smallest eigenvalue of ${\bf A}$.
\begin{thm}~\\
let ${\bf A} \stackrel{d}{=} KW_{p}(\nu, \Sigma)$, with $\nu=n-1 \geq p$, and let ${\bf \Lambda}>0$ be a $p\times p$ matrix. If $m=\frac{\nu-p-1}{2}$ is a positive integer, Then
\begin{eqnarray}\label{sarr1}
P({\bf A} > {\bf \Lambda})                 & = &  \frac{\Gamma\left(\frac{p(\nu+1)}{2}\right){\theta}^{\frac{p\nu}{2}}}{\Gamma\left(\frac{2q+p(\nu+1)-2}{2}\right)}\sum_{k=0}^{pm}{\sum_{\kappa}}^{\star}b_{k}\frac{( tr {\bf \Sigma}^{-1}{\bf \Lambda})^{\frac{p \nu}{2}-k}}{k!}C_{\kappa}({\bf
 \Sigma}^{-1}{\bf \Lambda})\nonumber \\
                                           & = & \frac{\Gamma\left(\frac{p(\nu+1)}{2}\right)}{\Gamma\left(\frac{2q+p(\nu+1)-2}{2}\right)}\sum_{k=0}^{pm}{\sum_{\kappa}}^{\star}b_{k}\frac{(\theta tr {\bf \Sigma}^{-1}{\bf \Lambda})^{\frac{p \nu}{2}-k}}{k!}C_{\kappa}(\theta{\bf\Sigma}^{-1}{\bf \Lambda})
\end{eqnarray}
where $\sum_{\kappa}^{\star}$ denotes summation over the partitions $\kappa=(k_{1}, \ldots, k_{p})$ of k with $k_{1}\leq m $, and
$$
b_{k}={\bf
G}^{30}_{23}\left({\theta}tr{\bf \Sigma}^{-1}{\bf
\Lambda}\begin{tabular}{|ccc}
                           $ 0$,                & $\frac{p}{2}$      &  \\
                           $-\frac{p\nu}{2}+k$, & $q-1+\frac{p}{2}$, & 0 \\
                     \end{tabular} \right)~,
$$
where $G^{ts}_{qr}$ denotes the Meijer's ${\bf G}$-function.
\end{thm}
{\bf Proof:}~\\
The pdf of ${\bf A}$ (\ref{simpdf}) is $f({\bf A})=C_{1}|{\bf \Sigma}|^{-\nu/2}|{\bf A}| ^{\frac{\nu-p-1}{2}}h_{1}(tr{\bf \Sigma}^{-1}{\bf
A})$, where $h_{1}(y)=(\theta y)^{\frac{2q+p-4}{4}}\exp(-\frac{\theta y}{2})W_{\alpha,\beta}(\theta y)$. Thus
$$
P({\bf A} > {\bf \Lambda})=C_{1}|{\bf \Sigma}|^{-\nu/2}\int_{{\bf A} > {\bf \Lambda}}|{\bf A}| ^{\frac{\nu-p-1}{2}}h_{1}(tr{\bf \Sigma}^{-1}{\bf
A})(d{\bf A})
$$
Making the change of variables ${\bf X}={\bf A}-{\bf \Lambda}$ and putting $tr{\bf \Sigma}^{-1}{\bf
\Lambda}=d$ give:
$$
P({\bf A} > {\bf \Lambda})=C_{1}|{\bf \Sigma}|^{-\nu/2}\int_{{\bf X} > 0}|{\bf X}+{\bf \Lambda}| ^{\frac{\nu-p-1}{2}}h_{1}(tr{\bf \Sigma}^{-1}{\bf X}+d)(d{\bf X})
$$
On the other hand, considering (\ref{F}) under the conditions of Remark \ref{Fpoly}, with $m=\frac{\nu-p-1}{2}$ being a positive integer, enable writing
\begin{eqnarray*}
 |{\bf X}+{\bf \Lambda}|^{m} & = & |{\bf X}|^{m}{\cdot}_{1}F_{0}(-m;-{\bf X}^{-1}{\bf \Lambda})\\
                             & = &|{\bf X}|^{m} \sum_{k=0}^{pm}{\sum_{\kappa}}^{\star}(-1)^{k}(-m)_{\kappa}\frac{C_{\kappa}({\bf X}^{-1}{\bf \Lambda})}{k!},
\end{eqnarray*}
where ${\sum_{\kappa}}^{\star}$ denotes summation over those partitions $\kappa=(k_{1},\ldots,k_{p})$ of k such that $k_{1} \leq m$.
Then, the probability $P({\bf A} > {\bf \Lambda})$ becomes:
\begin{eqnarray}\label{sarr2}
P({\bf A} > {\bf \Lambda}) & = & C_{1}|{\bf \Sigma}|^{-\nu/2}\sum_{k=0}^{pm}{\sum_{\kappa}}^{\star}\frac{(-1)^{k}(-m)_{\kappa}}{k!}\times \nonumber \\
                           &   &\int_{{\bf X} > 0}|{\bf X}| ^{m}h_{1}(tr{\bf \Sigma}^{-1}{\bf X}+d)C_{\kappa}({\bf X}^{-1}{\bf \Lambda})(d{\bf X}).
\end{eqnarray}
The last integral in (\ref{sarr2}), called {\bf J}, is evaluated using Li's result (\ref{LiR}), with $a=\frac{\nu}{2}$, $(-a+\frac{p+1}{2})_{\kappa}=(-m)_{\kappa}$, ${\bf Z}={\bf \Sigma}^{-1}$ and ${\bf U}={\bf \Lambda}$. That is,
\begin{eqnarray}\label{sarr3}
{\bf J} & = &  \int_{{\bf X} > 0}|{\bf X}| ^{m}h(tr{\bf \Sigma}^{-1}{\bf X})C_{\kappa}({\bf X}^{-1}{\bf \Lambda})(d{\bf X})\nonumber \\
        & = & \frac{(-1)^{k}\Gamma_{p}(\frac{\nu}{2})}{\Gamma(\frac{\nu p}{2}-k)(-m)_{\kappa}}{\cdot}{\gamma}{\cdot} |{\bf \Sigma}|^{\frac{\nu}{2}}C_{\kappa}({\bf \Lambda}{\bf \Sigma}^{-1}).
\end{eqnarray}
Finally, the desired result is obtained by substituting (\ref{mar1}) and (\ref{sarr3}) in (\ref{sarr2}), and using the fact that ${\theta}^{-k}C_{\kappa}(\theta {\bf \Sigma}^{-1}{\bf \Lambda}) = C_{\kappa}({\bf \Sigma}^{-1}{\bf \Lambda})$.
\begin{rmq}~\\
As a special case, if $q=2 \theta =1$, then $\alpha+\beta=\frac{1}{2}$ and using the identity (\ref{simpW}), we get a simple expression for $\gamma$:
\begin{eqnarray}
\gamma & = & {\theta}^{-\frac{\nu p}{2}+k}e^{-\theta d}\int_{0}^{\infty}t^{\frac{\nu p}{2}-k-1}e^{-\frac{t}{2}}dt\nonumber \\
       & = & (1/2)^{-\frac{\nu p}{2}+k}etr\left (-\frac{1}{2}{\bf \Sigma}^{-1}{\bf \Lambda}\right )\Gamma\left(\frac{\nu p}{2}-k\right).
\end{eqnarray}
Consequently, the expression (\ref{sarr2}) reduces to:
\begin{eqnarray}\label{mar2}
P({\bf A} > {\bf \Lambda}) & = & \sum_{k=0}^{pm}{\sum_{\kappa}}^{\star}\frac{(1/2)^{k}e^{-d/2}}{k!}C_{\kappa}({\bf \Sigma}^{-1}{\bf \Lambda})\nonumber \\
                           & = & etr\left(-\frac{1}{2}{\bf \Sigma}^{-1}{\bf \Lambda}\right)\sum_{k=0}^{pm}{\sum_{\kappa}}^{\star}\frac{C_{\kappa}\left(\frac{1}{2}{\bf \Sigma}^{-1}{\bf \Lambda}\right)}{k!}~,
\end{eqnarray}
which is similar to Khatri's result (\cite{khatri}). It is worth noting here that Khatri (\cite{khatri}) established his result for the random matrix ${\bf S} = \frac{1}{2}{\bf A}$, where ${\bf A}~\sim~W_{p}(\nu,\Sigma)$.
\end{rmq}
Two direct consequences of the above Theorem are provided below.
\begin{cor}~\\
let ${\bf A} \stackrel{d}{=} KW_{p}(\nu, \Sigma)$, with $\nu=n-1 \geq p$, and let ${\omega}_{p}$ be the smallest eigenvalue of ${\bf A}$. {\bf If} $m=\frac{\nu-p-1}{2}$ is a positive integer, Then
\begin{equation}\label{sarr4}
P({\omega}_{p} > x) =
\frac{\Gamma\left(\frac{p(\nu+1)}{2}\right)}{\Gamma\left(\frac{2q+p(\nu+1)-2}{2}\right)}x^{\frac{p\nu}{2}}\sum_{k=0}^{pm}{\sum_{\kappa}}^{\star}\frac{b_{k}}{k!}(\theta tr{\bf \Sigma}^{-1})^{\frac{p \nu}{2}-k}C_{\kappa}(\theta{\bf \Sigma}^{-1})
\end{equation}
where $\sum_{\kappa}^{\star}$ denotes summation over the partitions $\kappa=(k_{1}, \ldots, k_{p})$ of k with $k_{1}\leq m $, and
$$
b_{k}={\bf
G}^{30}_{23}\left({\theta}x tr{\bf \Sigma}^{-1} \begin{tabular}{|ccc}
                                                $ 0$,                & $\frac{p}{2}$      &  \\
                                                $-\frac{p\nu}{2}+k$, & $q-1+\frac{p}{2}$, & 0 \\
                                                \end{tabular} \right)~,
$$
where $G^{ts}_{qr}$ denotes the Meijer's ${\bf G}$-function. Thus, the distribution function of ${\omega}_{p}$ is given by
\begin{equation}\label{sarr5}
F_{{\omega}_{p}}(x) = 1-\frac{\Gamma\left(\frac{p(\nu+1)}{2}\right)}{\Gamma\left(\frac{2q+p(\nu+1)-2}{2}\right)}x^{\frac{p \nu}{2}}\sum_{k=0}^{pm}{\sum_{\kappa}}^{\star}\frac{b_{k}}{k!}(\theta tr {\bf \Sigma}^{-1})^{\frac{p \nu}{2}-k}C_{\kappa}(\theta{\bf\Sigma}^{-1}).
\end{equation}
\end{cor}
{\bf Proof:}~\\ Since the inequality ${\omega}_{p}> x$ is equivalent to ${\bf A} > x{\bf I}_{p}$, then (\ref{sarr4}) is obtained directly from (\ref{sarr3}) by putting ${\bf \Lambda} = x{\bf I}_{p}$, with $x>0$.
\begin{cor}~\\
Let ${\eta}_{1}$ be the largest eigenvalue of the inverted Kotz-Wishart matrix ${\bf A}^{-1}$. Its cdf $F_{{\eta}_{1}}(\cdot)$ is given by:
\begin{equation}\label{sarr6}
F_{{\eta}_{1}}(y)=P({\eta}_{1} \leq y)=\frac{\Gamma\left(\frac{p(\nu+1)}{2}\right)}{\Gamma\left(\frac{2q+p(\nu+1)-2}{2}\right)}y^{-\frac{p\nu}{2}}\sum_{k=0}^{pm}{\sum_{\kappa}}^{\star}\frac{b_{k}^{\prime}}{k!}(\theta tr{\bf \Sigma}^{-1})^{\frac{p \nu}{2}-k}C_{\kappa}(\theta{\bf \Sigma}^{-1}),
\end{equation}
where
$$
b_{k}^{\prime}={\bf
G}^{30}_{23}\left(\frac{\theta}{y} tr{\bf \Sigma}^{-1} \begin{tabular}{|ccc}
                                                $ 0$,                & $\frac{p}{2}$      &  \\
                                                $-\frac{p\nu}{2}+k$, & $q-1+\frac{p}{2}$, & 0 \\
                                                \end{tabular} \right)~,
$$
\end{cor}
{\bf Proof:}~\\
Since ${\omega}_{p}=\frac{1}{{\eta}_{1}}$, the expression (\ref{sarr6}) follows directly from (\ref{sarr4}) by putting $x=\frac{1}{y}$.\\
\subsubsection{\bf Some Comments}
\begin{itemize}
\item[(i)] The probability density function of the smallest eigenvalue ${\omega}_{p}$ of the KW matrix ${\bf A}$ may be derived by differentiating (\ref{sarr5}) with respect to $x$. It seems however difficult to express it in closed form.
\item[(ii)] Constantine \cite{const} established a closed form for the distribution function of the Wishart matrix, and then used it to derive explicit expression of the distribution function of the largest eigenvalue. In the next section, we will derive the distribution function of the inverted Kotz-Wishart matrix, which includes the cdf of the inverted Wishart matrix as a special case.
\end{itemize}
\subsection{\bf Distribution function of the Inverted Kotz-Wishart matrix}~\\
The distribution function of ${\bf A}^{-1}$ is obtained here using the following inequality in Loewner's sense (see \cite{haya}):
$$
{\bf A} > {\bf \Lambda}~is~equivalent~to~{\bf \Lambda}^{-1}> {\bf A}^{-1}> 0,
$$
this implies
\begin{equation}\label{Loewner}
P\left({\bf A} > {\bf \Lambda}\right)=P\left({\bf A}^{-1} < {\bf \Lambda}^{-1} \right)
\end{equation}
In other words, the cdf of ${\bf A}^{-1}$ may be derived from the expression $P\left({\bf A} > {\bf \Lambda}\right)$ that we have already evaluated.
\begin{thm}~\\
Let ${\bf B}={\bf A}^{-1}$ where ${\bf A} \stackrel{d}{=} KW_{p}(\nu, \Sigma)$, with $\nu=n-1 \geq p$, and let ${\bf \Omega}>0$ be a $p\times p$ matrix. {\bf If} $m=\frac{\nu-p-1}{2}$ is a positive integer, then the cdf of ${\bf B}$ is given by:
\begin{eqnarray}\label{sarr7}
P({\bf B} < {\bf \Omega})                  & = &  \frac{\Gamma\left(\frac{p(\nu+1)}{2}\right){\theta}^{\frac{p\nu}{2}}}{\Gamma\left(\frac{2q+p(\nu+1)-2}{2}\right)}\sum_{k=0}^{pm}{\sum_{\kappa}}^{\star}b_{k}^{\star}\frac{( tr {\bf \Sigma}^{-1}{\bf \Omega}^{-1})^{\frac{p \nu}{2}-k}}{k!}C_{\kappa}({\bf
 \Sigma}^{-1}{\bf \Omega}^{-1})\nonumber \\
                                           & = & \frac{\Gamma\left(\frac{p(\nu+1)}{2}\right)}{\Gamma\left(\frac{2q+p(\nu+1)-2}{2}\right)}\sum_{k=0}^{pm}{\sum_{\kappa}}^{\star}b_{k}^{\star}\frac{(\theta tr {\bf \Sigma}^{-1}{\bf \Omega}^{-1})^{\frac{p \nu}{2}-k}}{k!}C_{\kappa}(\theta{\bf\Sigma}^{-1}{\bf \Omega}^{-1})
\end{eqnarray}
where $\sum_{\kappa}^{\star}$ denotes summation over the partitions $\kappa=(k_{1}, \ldots, k_{p})$ of k with $k_{1}\leq m $, and
$$
b_{k}^{\star}={\bf
G}^{30}_{23}\left({\theta}tr{\bf \Sigma}^{-1}{\bf
\Omega}^{-1}\begin{tabular}{|ccc}
                           $ 0$,                & $\frac{p}{2}$      &  \\
                           $-\frac{p\nu}{2}+k$, & $q-1+\frac{p}{2}$, & 0 \\
                     \end{tabular} \right)~,
$$
\end{thm}
{\bf Proof:}~\\
The above result is proved by putting ${\bf \Lambda}={\bf \Omega}^{-1}$ in (\ref{sarr1}). Now, the cdf of the inverted Wishart matrix is obtained as a special case of (\ref{sarr7}).
\begin{cor}~\\
If $q=2\theta=1$ and $m=\frac{\nu-p-1}{2}$ is a positive integer, then the cdf of the inverted Wishart matrix is given by:
\begin{equation}\label{mar3}
P({\bf B} < {\bf \Omega})=
etr\left(-\frac{1}{2}{\bf \Sigma}^{-1}{\bf \Omega}^{-1}\right)\sum_{k=0}^{pm}{\sum_{\kappa}}^{\star}\frac{C_{\kappa}\left(\frac{1}{2}{\bf \Sigma}^{-1}{\bf \Omega}^{-1}\right)}{k!}.
\end{equation}
\end{cor}
{\bf Proof:} This is directly proved by putting ${\bf \Lambda}={\bf \Omega}^{-1}$ in (\ref{mar2}). We propose now an expression for the cdf $P({\bf A} < {\bf \Omega})$ of the Kotz-Wishart matrix {\bf A} as a {\bf conjecture:}\\ Let ${\bf A} \stackrel{d}{=} KW_{p}(\nu, \Sigma)$, with $\nu=n-1 \geq p$, and let ${\bf \Omega}>0$ be a $p\times p$ matrix. Then the cdf of {\bf A} is of the form:
\begin{eqnarray}\label{conj1}
P({\bf A} < {\bf \Omega}) & = & \int_{0}^{\Omega}f(A)(dA) \nonumber \\
                          & = & \frac{|{\theta}{\bf \Sigma}^{-1}{\bf \Omega}|^{\nu /2}{\Gamma}_{p}\left(\frac{p+1}{2}\right)}{{\Gamma}_{p}\left(\frac{\nu+p+1}{2}\right)}\sum_{k=0}^{\infty}\sum_{\kappa}\frac{(a)_{\kappa}}{(b)_{\kappa}}{\omega}_{k}(q,\theta)\frac{C_{\kappa}(-\theta{\bf \Sigma}^{-1}{\bf \Omega})}{k!},
\end{eqnarray}
where $a=\frac{\nu}{2}$, $b=\frac{\nu+p+1}{2}$ and ${\omega}_{k}(q,\theta)$ is a function satisfying the condition
\begin{equation}\label{conj2}
{\omega}_{k}(q=1,\theta=1/2)~=~1.
\end{equation}
When $q=2\theta=1$, the expression (\ref{conj1}) reduces to the cdf of the Wishart matrix, which was derived by Constantine \cite{const}. This is one of the reasons that motivated the conjecture.
\subsection{\bf The moment generating function of Kotz-Wishart matrix}~\\
Let ${\bf A} \stackrel{d}{=} KW_{p}(\nu, \Sigma)$, with $\nu=n-1 \geq p$, and let ${\bf \Omega}>0$ be a $p\times p$ matrix. The moment generating function $M_{\bf A}(\cdot)$ of ${\bf A}$ is defined by:
\begin{equation}
M_{\bf A}({\bf \Omega}) \stackrel{def}{=} E\left[etr({\bf \Omega}{\bf A})\right].
\end{equation}
Now using (\ref{F00}) and the linearity of the expectation-operator, we have
\begin{equation}\label{mgf}
M_{\bf A}({\bf \Omega})=\sum_{k=0}^{\infty}\sum_{\kappa}\frac{1}{k!}E\left(C_{\kappa}({\bf \Omega}{\bf A})\right)
\end{equation}
Before deriving the mgf of ${\bf A}$, let us first evaluate the expected value that appears in (\ref{mgf}).
\begin{lemma}~\\
Let ${\bf A} \stackrel{d}{=} KW_{p}(\nu, \Sigma)$, with $\nu=n-1 \geq p$, and let ${\bf \Omega}>0$ be a $p\times p$ matrix. And let k be a given positive integer, with $\kappa=(k_{1},\ldots,k_{p})$ a partition associated to k. Then we have:
\begin{equation}\label{expzon}
E\left[C_{\kappa}({\bf \Omega}{\bf A})\right]=K(n,p)C_{\kappa}({\bf \Omega}{\bf \Sigma}),
\end{equation}
where
\begin{equation}\label{coef1}
K(n,p)={\theta}^{-k}\frac{(\frac{n-1}{2})_{\kappa}\Gamma(\frac{np}{2})\Gamma(\frac{2q+np+2k-2}{2})}{\Gamma(\frac{2q+np-2}{2})\Gamma(\frac{np+2k}{2})}
\end{equation}
\end{lemma}
{\bf Proof:}~\\
Since the pdf of ${\bf A}$ is $f({\bf A})=C_{1}|{\bf \Sigma}|^{-\nu/2}|{\bf A}| ^{\frac{\nu-p-1}{2}}h_{1}(tr{\bf \Sigma}^{-1}{\bf
A})$, where $h_{1}(y)=(\theta y)^{\frac{2q+p-4}{4}}\exp(-\frac{\theta y}{2})W_{\alpha,\beta}(\theta y)$, we have:
\begin{eqnarray*}
E\left[C_{\kappa}({\bf \Omega}{\bf A})\right] & = & \int_{\bf A>0}C_{\kappa}({\bf \Omega}{\bf A})f({\bf A})(d{\bf A})\\
                                              & = & C_{1}|{\bf \Sigma}|^{-\nu/2}\int_{{\bf A} > 0}C_{\kappa}({\bf \Omega}{\bf A})|{\bf A}| ^{\frac{\nu-p-1}{2}}h_{1}(tr{\bf \Sigma}^{-1}{\bf A})(d{\bf A}).
\end{eqnarray*}
The latter integral is evaluated using the identity (\ref{garcia}), that is
$$
\int_{{\bf A} > 0}C_{\kappa}({\bf \Omega}{\bf A})|{\bf A}| ^{\frac{\nu-p-1}{2}}h_{1}(tr{\bf \Sigma}^{-1}{\bf A})(d{\bf A})=\frac{(a)_{\kappa} \Gamma_{p}(a)}{\Gamma(ap+k)}{\cdot}{\delta}{\cdot} |{\bf Z}|^{-a}C_{\kappa}({\bf
U}{\bf Z}^{-1}),
$$
where $a=\frac{\nu}{2}$, ${\bf Z}={\bf \Sigma}^{-1}$ and ${\bf U}={\bf \Omega}$.
Consequently the proof is completed by calculating the constant $\delta$, which is done as follows:
\begin{eqnarray*}
\delta & = & \int_{0}^{\infty}(\theta y)^{\frac{2q+p-4}{4}}\exp(-\frac{\theta y}{2})W_{\alpha,\beta}(\theta y)y^{\frac{\nu p}{2}+k-1}dy \\
       & = &{\theta}^{-\frac{\nu p}{2}-k+1}\int_{0}^{\infty}(\theta y)^{\frac{2q+p-4}{4}+\frac{\nu p}{2}+k-1}\exp(-\frac{\theta y}{2})W_{\alpha,\beta}(\theta y)dy \\
       & = &{\theta}^{-\frac{\nu p}{2}-k}\int_{0}^{\infty}z^{\frac{2q+p-4}{4}+\frac{\nu p}{2}+k-1}\exp(-\frac{z}{2})W_{\alpha,\beta}(z)dz \\
       & = &{\theta}^{-\frac{\nu p}{2}-k}\frac{\Gamma(\varepsilon+\frac{1}{2}-\beta)\Gamma(\varepsilon+\frac{1}{2}+\beta)}{\Gamma(\varepsilon-\alpha+1)},
\end{eqnarray*}
where the last expression is obtained after the change of variables $ z={\theta}y$, followed by the use of the identity (\ref{Whit}). Here, the appropriate parameters are: $ \varepsilon= \frac{2q+p-4}{4}+\frac{\nu p}{2}+k$, $\varepsilon+\frac{1}{2}-\beta =\frac{np-p+2k}{2}$, $\varepsilon+\frac{1}{2}+\beta =\frac{2q+np+2k-2}{2}$ and $\varepsilon-\alpha +1 =\frac{np+2k}{2}$. Hence we obtain
\begin{equation}\label{delta}
\delta = {\theta}^{-\frac{\nu p}{2}-k}\frac{\Gamma(\frac{2q+np+2k-2}{2})\Gamma(\frac{np-p+2k}{2})}{\Gamma(\frac{np+2k}{2})}
\end{equation}
Finally, the desired result is obtained by substituting (\ref{delta}) in the expression of $E\left[C_{\kappa}({\bf \Omega}{\bf A})\right]$, followed by some trivial simplifications. \\Next we establish a closed form of the moment generating function of the Kotz-Wishart matrix.
\begin{thm}~\\
Let ${\bf A} \stackrel{d}{=} KW_{p}(\nu, \Sigma)$, with $\nu=n-1 \geq p$, and let ${\bf \Omega}>0$ be a $p\times p$ matrix. And let k be a given positive integer, with $\kappa=(k_{1},\ldots,k_{p})$ a partition associated with k. Then, the mgf of ${\bf A}$ is given by:
\begin{equation}\label{Amgf}
M_{\bf A}({\bf \Omega})=\sum_{k=0}^{\infty}\sum_{\kappa}K(n,p)\frac{C_{\kappa}({\bf \Omega}{\bf \Sigma})}{k!}~,
\end{equation}
where $K(n,p)$ is given in (\ref{coef1}).
\end{thm}
{\bf Proof:}~\\To obtain the desired result (\ref{Amgf}), it suffices to substitute (\ref{coef1}) and (\ref{expzon}) in (\ref{mgf}).\\
As a special case of (\ref{Amgf}), if $q=2\theta=1$, then $K(n,p)$ and $E\left[C_{\kappa}({\bf \Omega}{\bf A})\right]$ reduce respectively to:
$$
K(n,p)=2^{k}\left(\frac{n-1}{2}\right)_{\kappa},~~~and~~~E\left[C_{\kappa}({\bf \Omega}{\bf A})\right]= \left(\frac{n-1}{2}\right)_{\kappa}C_{\kappa}(2{\bf \Omega}{\bf \Sigma}).
$$
And the mgf (\ref{Amgf}) becomes:
\begin{eqnarray*}
M_{\bf A}({\bf \Omega}) & = & \sum_{k=0}^{\infty}\sum_{\kappa}\left(\frac{n-1}{2}\right)_{\kappa}\frac{C_{\kappa}(2{\bf \Omega}{\bf \Sigma})}{k!}\\
                        & = & _{1}F_{0}\left(\frac{n-1}{2}; 2{\bf \Omega}{\bf \Sigma}\right)\\
                        & = & |I_{p}-2{\bf \Omega}{\bf \Sigma}|^{-\frac{n-1}{2}},~~~for~~0<2{\bf \Omega}{\bf \Sigma} < {\bf I}_{p},
\end{eqnarray*}
the last line being the well-known mgf of Wishart matrix $W_{p}(n-1, {\bf \Sigma})$.
\section{Matrix variate Varma transform}~\\
In 1951, Varma \cite{varma} introduced a generalization of the well-known Laplace transform, namely
\begin{equation}\label{univarma}
g(t)=\int_{0}^{\infty}\left[(tx)^{\beta-\frac{1}{2}}e^{-\frac{tx}{2}}W_{\alpha,\beta}(tx)\right]\phi(x)dx,~~~t>0
\end{equation}
where $ W_{\alpha,\beta}(tx)$, the Whittaker functions are defined in (\ref{w2}). In particular, if $\alpha+\beta=\frac{1}{2}$, then (\ref{univarma}) reduces to the classical Laplace transform, by virtue of the identity (\ref{simpW}),
\begin{equation}\label{unilap}
g(t)=\int_{0}^{\infty}e^{-tx}\phi(x)dx,~~~~t>0
\end{equation}
Formally, (\ref{univarma}) becomes an integral equation if $g(\cdot)$ is given and $\phi(\cdot)$ is to be found.
The purpose of this section is to:
\begin{itemize}
\item[(i)]   Introduce a matrix variate version of Varma transform (\ref{univarma}), that will be called M-Varma transform,
\item[(ii)]  Show that, as in the univariate case, M-Varma transform becomes the matrix variate Laplace transform when $\alpha+\beta=\frac{1}{2}$.
\item[(iii)] Evaluate M-Varma transforms for some specific functions of matrix argument.
\item[(iv)] Extend a confluent hypergeometric function using the M-Varma transform operator.
\end{itemize}
Let us first recall the definition of the Laplace transform (see \cite[p.252]{muir})
\begin{dfn}~\\
If $\phi({\bf X})$ is a function of the positive definite $p\times p$ matrix ${\bf X}$, the Laplace transform of $\phi({\bf X})$ is defined to be:
\begin{equation}\label{mlaplace}
L_{\phi}({\bf Z})= g({\bf Z})= \int_{{\bf X} > 0}etr(-{\bf X}{\bf Z})\phi({\bf X})(d{\bf X}),
\end{equation}
where ${\bf Z}={\bf U}+i{\bf V}$ is a complex symmetric matrix, ${\bf U}$ and ${\bf V}$ are real, and it is assumed that the integral is absolutely convergent in the right half-plane $Re({\bf Z})={\bf U}>{\bf U}_{0}$ for some positive definite ${\bf U}_{0}$.
\end{dfn}
As pointed out in Constantine \cite{const}, the Laplace transform satisfies the following convolution property:\\
If $g_{1}({\bf Z})$ and $g_{2}({\bf Z})$ are the Laplace transforms of $f_{1}({\bf S})$ and $f_{2}({\bf S})$, then $g_{1}({\bf Z})g_{2}({\bf Z})$ is
the laplace transform of the function
\begin{equation}\label{convol}
f({\bf R})=\int_{0}^{\bf R}f_{1}({\bf S})f_{2}({\bf R}-{\bf S})(d{\bf S}),~~~for~~0 < {\bf S} < {\bf R}.
\end{equation}
An illustrative example of the use of Laplace transform is provided through the pdf of Wishart matrix:
\begin{epl}~\\
Let ${\bf A} \stackrel{d}{=} W_{p}(n-1,{\bf \Sigma})$ whose pdf is given by:
$$
f({\bf A})~=~\frac{|{\bf A}|
^{\frac{n-1}{2}-\frac{p+1}{2}}}{2^{\frac{p(n-1)}{2}}{\Gamma}_{p}\left(\frac{n-1}{2}\right)|
{\bf \Sigma}| ^{\frac{n-1}{2}}}\exp\left(-\frac{1}{2}tr({\bf
\Sigma}^{-1}{\bf A})\right)
$$
Using the fact that $\int_{{\bf X}>0}f({\bf X})(d{\bf X})=1$, and putting ${\bf Z}=\frac{1}{2}{\bf \Sigma}^{-1}$ show that the Laplace transform of the function $\phi({\bf X})= |{\bf X}|^{\frac{n-1}{2}-\frac{p+1}{2}}$ is
\begin{equation}\label{lap1}
g({\bf Z})={\Gamma}_{p}\left(\frac{n-1}{2}\right)|{\bf Z}|^{-\frac{n-1}{2}}.
\end{equation}
\end{epl}
Inspired by the above example and the particular form of the pdf of the Kotz-Wishart matrix (\ref{simpdf}), we introduce now a generalization of (\ref{mlaplace}) as follows:
\begin{dfn}~\\
If $\phi({\bf X})$ is a function of the positive definite $p\times p$ matrix ${\bf X}$, the M-Varma transform of $\phi({\bf X})$ is defined to be:
\begin{equation}\label{mvarma}
V_{\phi}({\bf Z})= g({\bf Z})= \int_{{\bf X} > 0}\left\{\left[tr{\bf Z}{\bf X}\right]^{\frac{2q+p-4}{4}}etr\left(-\frac{1}{2}{\bf Z}{\bf X}\right)W_{\alpha,\beta}(tr{\bf Z}{\bf X})\right\}\phi({\bf X})(d{\bf X}),
\end{equation}
where $\alpha=\frac{2q-p}{4}$, $\beta=\frac{2q+p-2}{4}$, ${\bf Z}={\bf U}+i{\bf V}$ is a complex symmetric matrix, ${\bf U}$ and ${\bf V}$ are real, and it is assumed that the integral is absolutely convergent in the right half-plane $Re({\bf Z})={\bf U}>{\bf U}_{0}$ for some positive definite ${\bf U}_{0}$.
\end{dfn}
\begin{itemize}
\item[(i)] In the special case where $q=1$, then $\alpha$ and $\beta$ satisfy the condition: $\alpha+\beta=\frac{1}{2}$. Hence, using the identity
    (\ref{simpW}) shows that (\ref{mvarma}) becomes (\ref{mlaplace}).
\item[(ii)] Also, if $p=1$, then $\alpha=\beta=\frac{2q-1}{4}$, and $\frac{2q-3}{4}=\beta-\frac{1}{2}$. Consequently, the expression (\ref{mvarma}) reduces to
\begin{equation}\label{m1varma}
g(z)=\int_{0}^{\infty}\left[(zx)^{\beta-\frac{1}{2}}e^{-\frac{zx}{2}}W_{\alpha,\beta}(zx)\right]\phi(x)dx,~~~~z>0,
\end{equation}
which is clearly identical to the expression (\ref{univarma}).
\end{itemize}
To sum up, the M-Varma transform (\ref{mvarma}) appears to be a well-defined generalization of the original Varma transform and the classical matrix variate laplace transform (\ref{mlaplace}) as well. Next, we will evaluate some examples of M-Varma transforms for some functions of matrix argument.
\begin{lemma}~\\
Let $\phi_{1}({\bf X})= |{\bf X}|^{\frac{n-1}{2}-\frac{p+1}{2}}$, where $n \geq p+1$. Then the M-Varma transform of $\phi_{1}({\bf X})$ is given by:
\begin{equation}\label{examp1}
g_{1}({\bf Z})= \frac{\Gamma\left(\frac{2q+np-2}{2}\right){\Gamma}_{p}\left(\frac{n-1}{2}\right)}{\Gamma\left(\frac{np}{2}\right)}|{\bf Z}|^{-\frac{n-1}{2}}
\end{equation}
\end{lemma}
{\bf Proof:}~\\
Using the pdf (\ref{simpdf}) with ${\bf Z}={\theta}{\bf \Sigma}^{-1}$ gives
\begin{equation*}
\int_{{\bf X} > 0}\left\{\left[tr{\bf Z}{\bf X}\right]^{\frac{2q+p-4}{4}}etr\left(-\frac{1}{2}{\bf Z}{\bf X}\right)W_{\alpha,\beta}(tr{\bf Z}{\bf X})\right\}{\phi}_{1}({\bf X})(d{\bf X})=
\end{equation*}
\begin{equation*}
\frac{\Gamma\left(\frac{2q+np-2}{2}\right){\Gamma}_{p}\left(\frac{n-1}{2}\right)}{\Gamma\left(\frac{np}{2}\right)}|{\bf Z}|^{-\frac{n-1}{2}}.
\end{equation*}
As it may be noted, putting $q=2\theta=1$ in the last expression enables us to get back to (\ref{lap1}).
\begin{lemma}~\\
Let $\phi_{2}({\bf X})= |{\bf X}|^{\frac{n-1}{2}-\frac{p+1}{2}}C_{\kappa}({\bf X})$, where $n \geq p+1$ and let $\kappa=(k_{1},\ldots, k_{p})$ be a partition associated with the zonal polynomial $C_{\kappa}({\bf A})$ . Then the M-Varma transform of ${\phi}_{2}({\bf X})$ is given by:
\begin{equation}\label{examp2}
g_{2}({\bf Z})= \frac{\Gamma\left(\frac{2q+np+2k-2}{2}\right){\Gamma}_{p}\left(\frac{n-1}{2}\right)\left(\frac{n-1}{2}\right)_{\kappa}}{\Gamma\left(\frac{np+2k}{2}\right)}|{\bf Z}|^{-\frac{n-1}{2}}C_{\kappa}({\bf Z}^{-1})
\end{equation}
\end{lemma}
{\bf Proof:} This follows from (\ref{expzon}), with ${\bf \Omega}={\bf I}_{p}$ and ${\bf Z}={\theta}{\bf \Sigma}^{-1}$.\\ An immediate consequence of (\ref{examp2}) leads us to a more general result: Assume that
\begin{eqnarray*}
{\phi}_{3}({\bf X}) & = & |{\bf X}|^{\frac{n-1}{2}-\frac{p+1}{2}}{\cdot}_{m}F_{r}(a_{1},\ldots, a_{m};b_{1},\ldots,b_{r};{\bf X})\\
                    & = & |{\bf X}|^{\frac{n-1}{2}-\frac{p+1}{2}}\sum_{k=0}^{\infty}\sum_{\kappa}\frac{(a_{1})_{\kappa}\cdots(a_{m})_{\kappa}}{(b_{1})_{\kappa} \cdots (b_{r})_{\kappa}}\frac{C_{\kappa}({\bf X})}{k!}.
\end{eqnarray*}
Then, the M-Varma transform of ${\phi}_{3}({\bf X})$ is given by:
\begin{eqnarray}\label{varmagen}
g_{3}({\bf Z})&=&\int_{{\bf X} > 0}\left\{\left[tr{\bf Z}{\bf X}\right]^{\frac{2q+p-4}{4}}etr\left(-\frac{1}{2}{\bf Z}{\bf X}\right)W_{\alpha,\beta}(tr{\bf Z}{\bf X})\right\}{\phi}_{3}({\bf X})(d{\bf X})\nonumber \\
              &=& {\Gamma}_{p}\left(\frac{n-1}{2}\right)|{\bf Z}|^{-\frac{n-1}{2}}\sum_{k=0}^{\infty}\sum_{\kappa}\frac{(a_{1})_{\kappa}\cdots(a_{m})_{\kappa}\left(\frac{n-1}{2}\right)_{\kappa}}{(b_{1})_{\kappa} \cdots (b_{r})_{\kappa}}{\omega}_{k}\frac{C_{\kappa}({\bf \Sigma}^{-1})}{k!},
\end{eqnarray}
where $${\omega}_{k}=\frac{\Gamma\left(\frac{2q+np+2k-2}{2}\right)}{\Gamma\left(\frac{np+2k}{2}\right)}.$$ When $q=2\theta=1$, then ${\omega}_{k}=1$, and therefore, (\ref{varmagen}) becomes
$$
{\Gamma}_{p}\left(\frac{n-1}{2}\right)|{\bf Z}|^{-\frac{n-1}{2}}{\cdot}_{m+1}F_{r}(a_{1},\ldots, a_{m}, a;b_{1},\ldots,b_{r};{\bf Z}^{-1}),
$$
where $a=\frac{n-1}{2}$. This last expression is similar to that given in \cite[p.260]{muir}, due to Herz \cite{herz}. In other words, the definition of the hypergeometric functions of matrix argument $_{m}F_{r}(\cdot)$ by means of the Laplace transform remains robust under the M-Varma transform.\\ The next example of the use of the M-Varma transform involves the generalized Laguerre polynomial which is defined as follows (see \cite[p.282]{muir}):
\begin{dfn}~\\
The generalized Laguerre polynomial $L_{\kappa}^{\gamma}({\bf X})$ of an $p \times p$ symmetric matrix {\bf X} corresponding to the partition $\kappa$ of k is
\begin{equation}
L_{\kappa}^{\gamma}({\bf X})=(\gamma +t)_{\kappa}C_{\kappa}({\bf I}_{p})\sum_{s=0}^{k}\sum_{o}\binom{\kappa}{o}\frac{C_{o}(-{\bf X})}{(\gamma+t)_{o}C_{o}({\bf I}_{p})}~~~~\gamma >-1,
\end{equation}
where the inner summation is over all partitions $o$ of the integer s and $t=\frac{p+1}{2}$. Further, the generalized binomial coefficients $\binom{\kappa}{o}$ is defined through the following identity:
\begin{equation}\label{genbin}
\frac{C_{\kappa}({\bf I}_{p}+{\bf Y})}{C_{\kappa}({\bf I}_{p})}=\sum_{s=0}^{k}\sum_{o}\binom{\kappa}{o}\frac{C_{o}({\bf Y})}{C_{o}({\bf I}_{p})}
\end{equation}
\end{dfn}
The laplace transform of the function ${\phi}_{4}({\bf X})=|{\bf X}|^{\gamma}L_{\kappa}^{\gamma}({\bf X})$ is provided in \cite[p.282]{muir}. The following lemma gives a generalization of that result with the help of the M-Varma transform.
\begin{lemma}~\\
The M-Varma transform of ${\phi}_{4}({\bf X})=|{\bf X}|^{\gamma}L_{\kappa}^{\gamma}({\bf X})$ is given by:
\begin{eqnarray}\label{laguer}
g_{4}({\bf Z})& = &\int_{{\bf X} > 0}\left\{\left[tr{\bf Z}{\bf X}\right]^{\frac{2q+p-4}{4}}etr\left(-\frac{1}{2}{\bf Z}{\bf X}\right)W_{\alpha,\beta}(tr{\bf Z}{\bf X})\right\}{\phi}_{4}({\bf X})(d{\bf X})\nonumber \\
              & = & (\gamma +t)_{\kappa}{\Gamma}_{p}(\gamma+t)|{\bf Z}|^{-\gamma-t}C_{\kappa}({\bf I}_{p})\sum_{s=0}^{k}\sum_{o}\binom{\kappa}{o}{\omega}_{s}\frac{C_{o}(-{\bf Z}^{-1})}{C_{o}({\bf I}_{p})}
\end{eqnarray}
 where {\bf Z} is an $p\times p$ complex symmetric matrix with $Re({\bf Z})>0$, $t=\frac{p+1}{2}$, $\gamma >-1$ and $$ {\omega}_{s}=\frac{\Gamma\left(\frac{2q+np+2s-2}{2}\right)}{\Gamma\left(\frac{np+2k}{2}\right)}.$$
\end{lemma}
{\bf Proof:}~\\ Let $\varphi({\bf Z},{\bf X})=\left[tr{\bf Z}{\bf X}\right]^{\frac{2q+p-4}{4}}etr\left(-\frac{1}{2}{\bf Z}{\bf X}\right)W_{\alpha,\beta}(tr{\bf Z}{\bf X})$ be the kernel of the M-Varma transform. Then
\begin{equation*}
\int_{{\bf X} > 0}\varphi({\bf Z},{\bf X}){\phi}_{4}({\bf X})(d{\bf X})=
\end{equation*}
\begin{equation}\label{lag1}
(\gamma +t)_{\kappa}C_{\kappa}({\bf I}_{p})\sum_{s=0}^{k}\sum_{o}\frac{\binom{\kappa}{o}(-1)^{s}}{(\gamma +t)_{o}C_{o}({\bf I}_{p})}\int_{{\bf X} > 0}\varphi({\bf Z},{\bf X})|{\bf X}|^{\gamma+t-t}C_{o}({\bf X})(d{\bf X})
\end{equation}
The last integral is evaluated using (\ref{examp2}):
\begin{equation}\label{lag2}
\int_{{\bf X} > 0}\varphi({\bf Z},{\bf X})|{\bf X}|^{\gamma+t-t}C_{o}({\bf X})(d{\bf X})={\omega}_{s}(\gamma +t)_{o}{\Gamma}_{p}(\gamma+t)|{\bf Z}|^{-\gamma-t}C_{o}({\bf Z}^{-1}),
\end{equation}
where ${\omega}_{s}=\frac{\Gamma\left(\frac{2q+np+2s-2}{2}\right)}{\Gamma\left(\frac{np+2k}{2}\right)}$.
Therefore the desired result (\ref{laguer}) is obtained by substituting (\ref{lag2}) in (\ref{lag1}).\\ In particular, if $q=2\theta=1$, then using (\ref{genbin}) reduces (\ref{laguer}) to the Laplace transform of the function ${\phi}_{4}({\bf X})=|{\bf X}|^{\gamma}L_{\kappa}^{\gamma}({\bf X})$, provided in \cite[p.282]{muir}.\\We end up this section by extending the following confluent hypergeometric function of matrix argument, defined in \cite[p.472]{muir}:
\begin{equation}\label{cfgeo}
\psi(a,c;{\bf X})=\frac{1}{{\Gamma}_{p}(a)}\int_{{\bf Y}>0}etr(-{\bf Y}{\bf X})|{\bf Y}|^{a-p}|I+{\bf Y}|^{c-a-p}(d{\bf Y}),
\end{equation}
where $Re({\bf X})>0$ and $Re(a)>\frac{p-1}{2}$. In fact, (\ref{cfgeo}) may be viewed as the Laplace transform of the function
$$
\phi({\bf Y})= \frac{1}{{\Gamma}_{p}(a)}|{\bf Y}|^{a-p}|I+{\bf Y}|^{c-a-p}.
$$
Since the M-Varma transform generalizes the laplace transform, the above confluent hypergeometric function may be naturally extended as follows:
\begin{equation}\label{gsarr}
{\psi}^{q}(a,c;{\bf X})= \int_{{\bf Y} > 0}\left\{\left[tr{\bf X}{\bf Y}\right]^{\frac{2q+p-4}{4}}etr\left(-\frac{1}{2}{\bf X}{\bf Y}\right)W_{\alpha,\beta}(tr{\bf X}{\bf Y})\right\}\phi({\bf Y})(d{\bf Y}).
\end{equation}
Therefore, (\ref{gsarr}) represent the M-Varma transform for the function $\phi({\bf Y})$. Once again, by virtue of the identity (\ref{simpW}), (\ref{gsarr}) simplifies to (\ref{cfgeo}) when $q=1$.\\
 Let us now list some natural questions and comments regarding the M-Varma transform.
\begin{itemize}
\item[(i)] The M-Varma transform has been introduced as a generalization of the matrix variate Laplace transform (\ref{mlaplace}), which admits an inversion formula. So, a natural question is: Does an inversion formula for M-Varma transform exist in a closed form?
\item[(ii)] If so, does it simplify to the inversion Laplace transform for appropriately selected values of the parameters?
\item[(iii)] By analogy with the Laplace transform, does the M-Varma transform enjoy the convolution property (\ref{convol})?
\item[(iv)] For the univariate case ($p=1$), inversion formulae for the Varma transform (\ref{univarma}) have been derived under some conditions; see \cite{sri} where the author made use of the fractional integration theory, with the help of the Kober's operators.
\item[(v)] Hence, extending the Kober's operators to some functions of matrix argument might be a good starting point to derive the needed inversion formula.
\end{itemize}
\section{Concluding Remarks}~\\
In this paper, we have introduced a generalized Wishart distribution, namely, the Kotz-Wishart distribution, formed through the Kotz type distribution. Several results concerning both Kotz-Wishart matrix and Inverted Kotz-Wishart matrix are obtained. More precisely, exact distributions are established for all the statistics treated in our work. In particular, we extended Khatri's results \cite{khatri} by expressing the distribution function of the smallest eigenvalue of the Kotz-Wishart matrix in closed form, and under the same assumptions made by the above quoted author.\\
Nowadays, the Wishart distributions found an impressive number of applications in diverse areas, including Physics, Economics, Engineering, Biology and, of course, Multivariate statistical analysis. Therefore, the proposed Kotz-Wishart distributions investigated in our present work could be a very good alternative to Wishart distributions, especially for the situations where the normality assumption of the observations is arguably rejectable. Further, the presence of the additional parameters $q$, $\theta$ and $s$ in our findings provides more flexibility in comparison to the classical normal model. It is worth noting here that the use of exact distributions involving zonal polynomials is now less problematic in applied contexts, since Koev and Edelman (2006) \cite{koev} developped a very efficient algorithm for this issue.\\
Finally, given the key role played by the Laplace and inverse Laplace transforms in multivariate statistical theory, establishing an inversion formula for the M-Varma transform with the help of generalized Kober's operators would be a nice and promising bridge between Statistics and the theory of Integral equations. In the meantime, conjectures about the inversion formula for the M-Varma transform are, of course, welcome.
\appendix
\section{The Whittaker Functions}\label{app} The Whittaker function $W_{\kappa
,\mu}(z)$ is solution of the differential equation
\begin{equation}\label{w1}
\frac{d^{2}W}{d{z}^2}+\left(-\frac{1}{4}+\frac{\kappa}{z}+\frac{\frac{1}{4}-{\mu}^{2}}{z^2}\right)W~=~0.
\end{equation}
This equation is called the Whittaker's equation;
$W_{\kappa,\mu}(z)$ can be expressed as follows (see
\cite[p.1025]{ryzhik})
\begin{equation}\label{w2}
W_{\kappa,\mu}(z)~=~\frac{z^{\kappa}e^{-\frac{z}{2}}}{\Gamma\left(\mu-\kappa+\frac{1}{2}\right)}\int_{0}^{\infty}t^{\mu-\kappa-\frac{1}{2}}e^{-t}\left(1+\frac{t}{z}\right)^{\mu+\kappa-\frac{1}{2}}dt,
\end{equation}
\begin{equation*}
\left[Re(\mu-\kappa)>-\frac{1}{2},~~\mid arg(z)\mid <\pi \right].
\end{equation*}
In the special case where $\mu+\kappa~=~\frac{1}{2}$, the above
expression simplifies to:
\begin{eqnarray}\label{simpW}
W_{\kappa,\mu}(z) & = & \frac{z^{\kappa}e^{-\frac{z}{2}}}{\Gamma\left(\mu-\kappa+\frac{1}{2}\right)}\int_{0}^{\infty}t^{\mu-\kappa+\frac{1}{2}-1}e^{-t}dt \nonumber \\
                  & = & z^{\kappa}e^{-\frac{z}{2}}.
\end{eqnarray}
There are many definitions and integral representations of the
Whittaker functions $W_{\kappa,\mu}(z)$ in the literature, (see
\cite[p.505]{abra}, \cite[p.1024]{ryzhik}).
Here, we just collect some results for $W_{\kappa,\mu}(\cdot)$ which
are used in our work. The following integrals involving Whittaker
function are taken from \cite[p.26]{ober}, \cite[p.823,827]{ryzhik}.
\begin{equation}\label{mellin}
\int_{0}^{\infty}(b+x)^{\nu}e^{-ax}x^{y-1}dx=b^{\frac{y+\nu-1}{2}}a^{-\frac{y+\nu+1}{2}}e^{\frac{ab}{2}}\Gamma(y){\cdot}W_{\alpha,\beta}(ab),
\end{equation}
where $W_{\alpha,\beta}(\cdot)$ denotes the Whittaker
function, with
$$
\left\{
\begin{array}{ccc}
\alpha & = & \frac{\nu-y+1}{2} \\
\beta  & = & \frac{\nu+y}{2}
\end{array}
\right.,
$$
and
$$
\Gamma(y)~=~\int_{0}^{\infty}x^{y-1}e^{-x}dx~.
$$
\begin{equation}\label{Meijer}
\int_{0}^{\infty}x^{\rho -1}(c+x)^{-\sigma}e^{-\frac{x}{2}}W_{\alpha,\beta}(c+x)dx=\Gamma(\rho)c^{\rho}e^{\frac{c}{2}}
G^{30}_{23}\left(c~\begin{tabular}{|ccc}
                           $ 0$,      & $1-\alpha-\sigma$        &  \\
                           $-\gamma$, & $\frac{1}{2}+\mu-\sigma$,&$\frac{1}{2}-\mu-\sigma$ \\
                     \end{tabular} \right)~,
\end{equation}
where $|arg(c)|< \pi$ , $Re(\rho) > 0$.
\begin{equation}\label{Whit}
\int_{0}^{\infty}e^{-\frac{1}{2}x}x^{\mu-1}W_{\alpha,\beta}(x)dx=\frac{\Gamma\left(\mu+\frac{1}{2}-\beta\right)\Gamma\left(\mu+\frac{1}{2}+\beta\right)}{\Gamma\left(\mu-\alpha+1\right)},
\end{equation}
where $Re(\mu+\frac{1}{2}\pm \beta)>0$.

\address{DEPARTMENT OF MATHEMATICS AND STATISTICS\\
SULTAN QABOOS UNIVERSITY\\
AL-KHOUD 123, MUSCAT, THE SULTANATE OF OMAN\\
\printead{e1}\\
\phantom{E-mail:asarr7@gmail.com ;asarr@squ.edu.om }}
\end{document}